\documentclass[11pt,a4paper,twoside,listof=totoc,usenames,dvipsnames]{article}
\usepackage{amsthm}
\usepackage{framed}
\usepackage{multirow}
\usepackage{amssymb}
\usepackage[utf8]{inputenc}
\usepackage{url}
\usepackage{authblk}

\usepackage{mathtools}
\usepackage{etoolbox}
\usepackage{float}
\usepackage{xcolor}
\usepackage{relsize}

\usepackage{calc}

\newcommand{\beginsupplement}{%
	\setcounter{table}{0}
	\renewcommand{\thetable}{S\arabic{table}}%
	\setcounter{figure}{0}
	\renewcommand{\thefigure}{S\arabic{figure}}%
}

\usepackage[
per-mode=fraction,
locale=US,
decimalsymbol=.,
obeyall,
]{siunitx}
\usepackage[utf8]{inputenc}
\usepackage[T1]{fontenc}

\usepackage{algpseudocode} 
\usepackage{algorithm} 

\DeclareSIUnit[number-unit-product = \,]\cells{\text{cells}}
\DeclareSIUnit[number-unit-product = \,]\Cells{\num{e5}\,\text{cells}}
\DeclareSIUnit[number-unit-product = {}]\fbs{\%\,\text{FBS}}
\DeclareSIUnit[number-unit-product = \,]\FBS{\num{10}\%\,\text{FBS}}
\DeclareSIUnit[number-unit-product = \,]\day{\text{day}}

\usepackage{graphicx}
\usepackage{subfigure}
\usepackage{bbm}

\usepackage{tabularx}

\usepackage{extarrows}
\usepackage{bm}
\usepackage{nicefrac}
\usepackage{cancel}
\usepackage[normalem]{ulem}

\newcolumntype{C}[1]{>{\centering\arraybackslash}p{#1}}

\usepackage{booktabs}

\usepackage{blindtext}
\usepackage{caption}

\usepackage{rotating}

\usepackage[toc,page]{appendix}

\usepackage[nottoc]{tocbibind}

\usepackage{stackengine}
\stackMath

\numberwithin{equation}{section}
\def\Triang{\mbox{Triang}}
\def\totDiff#1{\frac{\text{d}}{\text{d}#1}}
\def\partDiff#1{\frac{\partial}{\partial#1}}

\def\iBG{I^\text{BG}}
\def\initS{S_0}
\def\initV{V_0}
\def\stress{\eta}
\def\initStress{\eta_0}

\def\Sthr{S_{\text{thr}}}
\def\dV{\dot{V}}
\def\dStress{\dot{\stress}}
\def\equibV{\bar{V}}

\def\nat{\lambda}
\def\starv{\nat_{\text{st}}}
\def\scaledDeath{\nat_S}
\def\scalDeathStress{\nat_\stress}

\def\prol{\beta}
\def\scaledProl{\prol_S}
\def\netProl{\tilde{\prol}}
\def\scalProlStress{\prol_\stress}

\def\Kv{K}
\def\netKv{\tilde{K}}

\def\adaptS{\alpha_S}

\def\deactFct#1{\delta ^{-}(#1)}
\def\actFct#1{\delta^{+}(#1)}

\def\odeVprol{\prol\,V\left(1-\left(\frac{V}{\Kv }\right)^{m}\right)-\nat\,V}

\def\solLogistic#1#2{\frac{\initV#2}{\left(\initV^m+\left(#2^m-\initV^m\right)\cdot e^{-#1mt}\right)^{\frac{1}{m}}}}
\def\odeVstarv{-(\nat+\starv)V}

\def\IR{\mathbb{R}}

\def\mean{\mathbb{E}}
\def\Var{\text{Var}}
\def\unif{\mathcal{U}}

\numberwithin{equation}{section}

\begin{document}

\title{Environmental stress level to model tumor\\ cell growth and survival}

\author[1]{Sabrina Schönfeld\thanks{sabrina.schoenfeld@tum.de}}
\author[2]{Alican Ozkan}
\author[3]{Laura Scarabosio}
\author[4]{Marissa Nichole Rylander}
\author[1]{Christina Kuttler}

\affil[1]{\small Center of Mathematics, Technical University of Munich, Garching, Germany}
\affil[2]{\small Wyss Institute for Biologically Inspired Engineering, Harvard University, Boston, MA, United States}
\affil[3]{\small Institute for Mathematics, Astrophysics and Particle Physics, Radboud University, Nijmegen, Netherlands}
\affil[4]{\small Department of Mechanical Engineering, The University of Texas, Austin, TX, United States}

\renewcommand\Authands{ and }
\date{}
\maketitle

\paragraph{Abstract.} Survival of living tumor cells underlies many influences such as nutrient saturation, oxygen level, drug concentrations or mechanical forces. Data-supported mathematical modeling can be a powerful tool to get a better understanding of cell behavior in different settings. However, under consideration of numerous environmental factors mathematical modeling can get challenging. We present an approach to model the separate influences of each environmental quantity on the cells in a collective manner by introducing the "environmental stress level". It is an immeasurable auxiliary variable, which quantifies to what extent viable cells would get in a stressed state, if exposed to certain conditions. A high stress level can inhibit cell growth, promote cell death and influence cell movement. As a proof of concept, we compare two systems of ordinary differential equations, which model tumor cell dynamics under various nutrient saturations respectively with and without considering an environmental stress level. Particle-based Bayesian inversion methods are used to quantify uncertainties and calibrate unknown model parameters with time resolved measurements of \textit{in vitro} populations of liver cancer cells. The calibration results of both models are compared and the quality of fit is quantified. While predictions of both models show good agreement with the data, there is indication that the model considering the stress level yields a better fitting. The proposed modeling approach offers a flexible and extendable framework for considering systems with additional environmental factors affecting the cell dynamics.

\section{Introduction}
As cancerous diseases still cause many deaths in human population, there is a high need of improving the understanding of tumor growth dynamics and treatment strategies. For that purpose and for treatment optimization, mathematical approaches can help with adequate quantitative descriptions and predictions. Numerous mathematical models and methods are already available (see~\cite{Preziosi, Byrne}), including differential equations of all types, stochastic models, phenomenological models as well as more mechanistic ones. Despite good new experimental techniques, e.g. measuring tumor growth also {\it in vivo}~\cite{meas1,meas2,meas3}, it is still difficult to quantify growth by determining parameter values in realistic situations without harming a patient or affecting a potential treatment. 
To get a better quantitative understanding of the underlying processes, investigating tumor cell cultures and even single cells under well-controlled conditions provides a promising way to get at least parts of the necessary information (see e.g.~\cite{Szot}). But even under laboratory conditions, parameters can underlie variations and cells may behave individually to some extent. To consider this, uncertainty quantification can be an adequate tool. Having such {\it in vitro} experiments, it is much easier to separate and measure different environmental influences on the tumor growth. The probably most relevant environmental influences are the availability of nutrients and oxygen. For modeling purposes, one often focuses on the limiting nutrient, without specifying it in detail. Nutrient deprivation causes ``stress'' to any cells, but other factors may cause stress to the cells, as well. Thus, it is worth to consider the concept of stress as a relevant factor for modeling growth and death of tumor cells.

To approach this problem, we suggest the introduction of an ``environmental stress level'' as an immeasurable auxiliary variable, which quantifies to what extent viable cells would get stressed out under certain environmental conditions. These may include quantities like nutrient saturation, oxygen level, drug concentrations or mechanical forces. On the one hand, this environmental stress level can e.g. inhibit cell growth and promote death by inducing necrosis for present cells. On the other hand, it is assumed to have an influence on cell movement, as cells favor areas with a low environmental stress level, which may cause metastasis formation, as well. In this work, this approach is presented in a simplified spatially homogeneous setting, where only the nutrient saturation effects the environmental stress level. Nevertheless, the presented model can be easily extended to more environmental variables, which provides several advantages from a modeling point of view. Even in this very simplified setting, model comparison shows a slight preference for the newly introduced model.

For comparison, we present two systems of ordinary differential equations, which model tumor cell population dynamics respectively with and without considering an environmental stress level. To do so, time-resolved \emph{in vitro} data (from~\cite{Lima}) are used to perform model parameter calibration with Bayesian inversion. The latter approach provides a quantitative estimate of the uncertainty in the estimated parameters and hence brings more information than a deterministic inversion method, at the price of a higher computational cost. In the Bayesian inversion setting several methods are available for sampling. Among the most popular ones are Markov Chain Monte Carlo (MCMC) methods~\cite[Ch. 6-7]{RC04}, which are very simple to implement, but they need ad hoc tuning when the dimension of the parameter space is high or for complex posterior distributions~\cite{propKernel}. In the last decade, particle based methods have gained popularity for their higher efficiency and robustness compared to simple MCMC, especially in the presence of high dimensional parameter spaces and multimodal posterior distributions~\cite{evid}. Most of these algorithms are population based, which means they deal with a collection of samples in every iteration. In this category are sequential importance sampling and resampling (annealed importance sampling~\cite{Neal01} and Sequential Monte Carlo \cite{C02,SMC}) and population MCMC~\cite{LW01,JSH07}. To calibrate the unknown parameters of the presented models, we use Sequential Monte Carlo (SMC). It is worth mentioning that the ensemble Kalman filter~\cite{E03} can also be used to solve the inverse problem. It can be viewed as an approximate SMC with Gaussian approximation of the posterior~\cite{SS17}. However, the rigorous analysis of this method is confined for the moment to linear forward operators~\cite{BSWW19,SS17}. Eventually, the calibration results are used to compare the models via a validation metric and the Bayes factor. For this purpose, the chosen SMC method is advantageous opposed to other algorithms (e.g. MCMC), since its structure provides easy access to the Bayes factor~\cite{evid}.

The article is structured as follows: The upcoming Section~\ref{sec:methods} explains the investigated mathematical models with the underlying experimental setup and provides background about the applied methods for uncertainty quantification, parameter calibration and model comparison. The results from the model calibrations are summarized and discussed in Section~\ref{sec:results}. Finally, Section~\ref{sec:conclusion} concludes the presented results and the use of the novel modeling approach. At the end of this article, more details about mathematical analysis of the models and the calibration results are enclosed in a supplementary section.

\section{Materials and method}
\label{sec:methods}
Section~\ref{subsec:model} presents different models for population dynamics of tumor cells and introduces the idea of using a stress level to model environmental influences on the cells collectively. To investigate the effect of nutrient changes on viable cells, we consider a special case of the models. The underlying biological setting and the measurement methods for the corresponding \textit{in vitro} experiments are explained in Section~\ref{subsec:experiments}. With these data we calibrate the unknown parameters of the models using Bayesian inversion and Sequential Monte Carlo methods. Section~\ref{subsec:calibration} provides the theoretical background behind the algorithms. Finally, we present the validation metric and the Bayes factor in Section~\ref{subsec:modelComp}, which are used to investigate and compare the resulting model calibrations.
\subsection{Mathematical modeling}
\label{subsec:model} 
The presented ODE models describe the basic growth and death dynamics of tumor cells in a spatially homogeneous, avascular environment under consideration of the present nutrient saturation. The time-dependent variables are the density of viable tumor cells~${V=V(t)}$, the nutrient saturation~${S=S(t)}$, and the environmental stress level~${\stress=\stress(t)}$ at day~${t\geq 0}$.

We consider two possibilities to model the influence of the nutrient saturation on the cells: Model~\ref{eq:ODE} assumes that the nutrient saturation directly affects the growth/death rates of the cells, while in model~\ref{eq:ODEstress} the nutrient supply influences the cell population indirectly by changing the environmental stress level, which itself has an effect on the viability of the cells. The nutrient saturation is bounded and normalized: ${S(t)\in[0,1]}$. The bounds represent the complete absence of nutrients (${S=0}$) and a nutrient supply, which creates optimal growth conditions for the cells ($S=1$). In the latter case, both models can be reduced to a simplified model, which is independent of the nutrient saturation~$S$. We denote this model by~\ref{eq:ODEprol} (``opt'' short for ``optimal'' growth conditions).

A detailed description of each model follows in the consecutive paragraphs. All models are biologically reasonable in a first check, as they have been mathematically analyzed in terms of positivity and boundedness of the solutions as well as steady states and their stability. The last paragraph of this section provides a summary of the variables, parameters, and mathematical features of the models.
\paragraph{What if the cells find optimal nutrient conditions?} External nutrients serve as an energy source for cell proliferation. Growth can happen with a maximal possible proliferation rate~$\prol$, which is virtually reached under optimal nutrient conditions (${S=1}$). Under these circumstances, cell death is assumed to occur only for nutrient-independent reasons (e.g. dying of old age) with a constant rate~$\nat$. These dynamics can be modeled with a combination of a generalized logistic growth term (first proposed in~\cite{genLog}) and an exponential death term. The population growth is limited by the carrying capacity~$\Kv$ of the biological system, taking into account e.g. limited space but not nutrient shortage. In particular, cell growth is inhibited in a cell density dependent manner, which is known as ``proliferation contact inhibition''~\cite{contactInhib}. The strength of this phenomenon is specific for the cell type and is modeled by the shape parameter~$m$ of the logistic growth. A small value for~$m$ indicates strong contact inhibition, i.e. reduction of proliferation already starts at small cell densities. Therefore, $1/m$ can be interpreted as the strength of contact inhibition. Overall, this results in the following initial value problem:
	\begin{align}
		&\!S=1~~\Rightarrow~~\left\{~~\begin{aligned}
			\dV 	&= \odeVprol\,, \\
			V(0) 	&= V_0\,.
		\end{aligned} \right. \tag{$\mathcal{M}_{\text{opt}}$} & \label{eq:ODEprol}
	\end{align}
As a feature of the pure logistic term, population growth for very small populations (i.e. ${V\ll\Kv}$) can be approximated by exponential growth with rate~$\prol$. It is a reasonable assumption that, independently of the initial population size~$\initV$, the cell population's size actually increases over time under optimal nutrient conditions, which translates to a parameter relation: ${\prol>\nat}$. With this constraint, the ODE can be rewritten as a purely logistic growth term with a new rate~$\netProl$ and capacity~$\netKv$:
	\begin{align*}
	\dV	&= \prol\,V\left(1-\frac{\nat}{\prol}-\left(\frac{V}{\Kv}\right)^m\right)
	\overset{\prol>\nat}{=} \prol\,V\left(1-\frac{\nat}{\prol}\right)\cdot\left(1-\left(\frac{V}{\Kv\left(1-\frac{\nat}{\prol}\right)^{\frac{1}{m}}}\right)^m\right) \\
	&= \underbrace{\left(\prol-\nat\right)}_{=\,\netProl}V\left(\vphantom{\frac{V}{\Kv\left(1-\frac{\nat}{\prol}\right)^{\frac{1}{m}}}}\right.1-\left(\vphantom{\frac{V}{\Kv\left(1-\frac{\nat}{\prol}\right)^{\frac{1}{m}}}}\right.\underbrace{\frac{V}{\Kv\left(1-\frac{\nat}{\prol}\right)^{\frac{1}{m}}}}_{=\,\netKv}\left.\vphantom{\frac{V}{\Kv\left(1-\frac{\nat}{\prol}\right)^{\frac{1}{m}}}}\right)^m\left.\vphantom{\frac{V}{\Kv\left(1-\frac{\nat}{\prol}\right)^{\frac{1}{m}}}}\right)\,.
	\end{align*}
Therefore, under optimal nutrient conditions, the cell population actually grows logistically with
	\begin{align*}
	\mbox{``net'' growth rate:}\quad \netProl&=(\prol-\lambda)\in(0,\prol) \\
	\mbox{and ``net'' carrying capacity:}\quad \netKv&=\Kv\left(1-\frac{\nat}{\prol}\right)^{\frac{1}{m}}=\Kv\left(\frac{\prol-\nat}{\prol}\right)^{\frac{1}{m}}\in(0,\Kv)\,.
	\end{align*}
For~$t\geq 0$ the known analytical solution of such an initial value problem, i.e. for model~\ref{eq:ODEprol}, is
\begin{align}
	V(t)&=\solLogistic{\netProl}{\netKv}\nonumber \\
	&=\left(\frac{(\initV\Kv)^m(\prol-\nat)}{{\prol}\initV^m+\left(\left(\prol-\nat\right)\Kv^m-\prol\initV^m\right)\cdot e^{-\left(\prol-\nat\right)mt}}\right)^{\frac{1}{m}}\,. \label{eq:solV}
\end{align}
\paragraph{What if suboptimal nutrient conditions directly affect the cell dynamics?} Exposing the cells to suboptimal nutrient conditions (${0\leq S<1}$) can affect the cells' metabolism by decelerating proliferation or even lead to cell death by starvation. To include this, an additional death term is needed, whose rate as well as the proliferation rate are scaled appropriately to the nutrient supply.

\begin{figure}[H]
	\begin{center}
		\includegraphics[width=.4\textwidth]{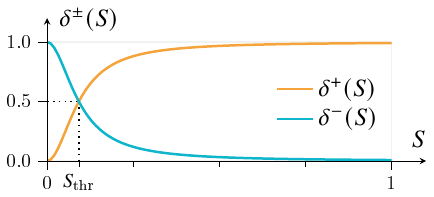}
		\caption{Qualitative behavior of the Hill type scaling functions~$\actFct{S}$ and~$\deactFct{S}$}
		\label{Fig:scalFct}
	\end{center}
\end{figure}

A time-dependent nutrient supply can be described by adding an additional differential equation to the system: ${\dot{S}=g(S,V,t)}$ with ${S(0)=\initS}$\,. The reaction term~$g(S,V,t)$ can e.g. include terms for nutrient consumption by viable cells or an external nutrient source. In the considered \emph{in vitro} experiments, whose measurements are used to calibrate the model later on, the nutrient concentration is maintained approximately constant. Therefore, all presented models assume 
	\begin{equation*}
	\dot{S}=g(S,V,t)=0 \quad\Rightarrow \quad S(t)=\mbox{const.}=\initS\in[0,1]\quad\forall t\geq 0 \,.
	\end{equation*}
The present nutrient saturation~$\initS$ contributes to the cells' reproductivity as well as to their ability to survive. If there are no nutrients available, the cells die from starvation with a maximal possible rate~$\starv$\,, assumed to happen in an exponential manner. With increasing nutrient supply, the starvation rate~$\starv$ decreases and the proliferation rate~$\prol$ is up-regulated. Adapting the previous model~\ref{eq:ODEprol} accordingly yields the more general initial value problem
	\begin{align}
	&\!\left\{~\begin{aligned}
		\dV &= \underbrace{\actFct{\initS}\cdot\prol}_{=\,\scaledProl}\, V\left(1-\left(\frac{V}{\Kv}\right)^{m}\right)-\underbrace{(\nat+\deactFct{\initS}\cdot\starv)}_{=\,\scaledDeath} V\,,\\[-4pt]
		V(0) &=\initV\,,
	\end{aligned} \right.\tag{$\mathcal{M}_S$} & \label{eq:ODE}
\end{align}
with appropriate nutrient-dependent functions~$\actFct{S}$ and~$\deactFct{S}$, scaling the total growth and death rates:
\begin{equation*}
	\scaledProl=\actFct{\initS}\cdot\prol\qquad\mbox{and}\qquad \scaledDeath=(\nat+\deactFct{\initS}\cdot\starv)\,.
\end{equation*} 
The scaling functions~$\delta^{\pm}(S)$ need to fulfill certain criteria to be biologically meaningful in our model setting. First, the cells should not die from starvation but proliferate with the maximal possible rate under optimal nutrient conditions.  In contrast, without nutrients the cells should not be able to reproduce but starve with a maximal possible rate. These features translate to the properties
\begin{equation*}
	\actFct{S}\cdot\prol~
	\begin{cases}
		~\to \prol &\mbox{for }S\to1\,,\\
		~\in (0,\prol) &\mbox{for }0<S<1\,,\\
		~\to 0 &\mbox{for }S\to0\,,
	\end{cases}
	\quad\mbox{and}\quad
	\deactFct{S}\cdot\starv~
	\begin{cases}
		~\to 0 &\mbox{for }S\to1\,,\\
		~\in (0,\starv) &\mbox{for }0<S<1\,,\\
		~\to \starv &\mbox{for }S\to0\,.
	\end{cases}
\end{equation*} 
Second, increasing the nutrient saturation should promote growth and reduce starvation: $\actFct{\initS}$ resp. $\deactFct{\initS}$ need to be monotonically increasing resp. decreasing functions bounded by~$[0,1]$. Therefore, a possible choice are Hill type functions with Hill coefficient~$k=2$, i.e.
	\begin{equation}
	\delta^+:S\mapsto\frac{S^2}{\Sthr^2+S^2} \qquad\mbox{and}\qquad \delta^-:S\mapsto 1-\actFct{S}=1-\frac{S^2}{\Sthr^2+S^2}\,, \label{eq:deltas}
	\end{equation}
where ${0<\Sthr\ll1}$ denotes the nutrient threshold, for which the cells proliferate/starve with half-maximal rate. For simplicity, we do not assume any hysteresis effect, motivating the relation~${\deactFct{S}=1-\actFct{S}}$. The scaling functions~$\delta^\pm(S)$ can be interpreted as influence functions, which describe how tolerant viable cells are to nutrient changes. In this context, the parameter~$\Sthr$ can be seen as a nutrient sensitivity threshold. Figure~\ref{Fig:scalFct} shows a qualitative plot of the behavior of the scaling functions.
It can be shown (for calculations see Section~\ref{SUPsec:analSol} in the supplement) that for~$t\geq 0$ a well-defined analytical solution of problem~\eqref{eq:ODE} is
\begin{equation*}
	V(t)
	= \begin{cases}
		~\initV\Kv\cdot\left(\frac{\scaledProl-\scaledDeath}{\scaledProl\initV^m+\left(\left(\scaledProl-\scaledDeath\right)\Kv^m-\scaledProl\initV^m\right)e^{-\left(\scaledProl-\scaledDeath\right)mt}}\right)^{\frac{1}{m}} &\mbox{if }\scaledProl\neq\scaledDeath\,, \\[8pt]
		~\initV\Kv\cdot\Big(\frac{1}{mt\scaledProl\initV^m+\Kv^m}\Big)^{\frac{1}{m}}
		&\mbox{if }\scaledProl=\scaledDeath\,.
	\end{cases}
\end{equation*}
We note that in the extreme situation of having no nutrient supply (i.e. ${\initS=0}$), it holds ${\actFct{\initS}=0}$ and  ${\deactFct{\initS}\approx 1}$. This results in ${\scaledProl=0}$ and ${\scaledDeath\approx\nat+\starv>0}$. By inserting this special case of~${\scaledProl\neq\scaledDeath}$ into the given analytical solution, it simplifies to
\begin{equation*}
	V(t)\overset{\initS\to 0}{\longrightarrow} \initV\Kv\cdot\left(\frac{-(\nat+\starv)}{-(\nat+\starv)\Kv^me^{(\nat+\starv)mt}}\right)^{\frac{1}{m}}=\initV e^{-(\nat+\starv)t}\,.
\end{equation*} 
This is an expected result, being the solution of a simple exponential decay model without growth:
\begin{align*}
	&\left\{~\begin{aligned}
		\dV &= \odeVstarv\,, \\
		V(0) &= \initV\,.
	\end{aligned} \right.
\end{align*} 
\paragraph{Advantages of using the environmental stress level (ESL) instead to directly affect the cells.} 
For an alternative way to include the effect of the nutrient conditions in model~\ref{eq:ODEprol}, we introduce an auxiliary variable: the environmental stress level~${\stress=\stress(t)}$. It is designed to be an immeasurable quantity, describing to what extent viable cells would get stressed out, if exposed to certain environmental conditions, like in this case the nutrient saturation. We assume that the stress level is limited by a maximal level, which is set to one: ${\stress\in[0,1]}$.

In a general setting, the the ESL can be influenced by multiple environmental factors (e.g nutrient/oxygen/drug concentration), where each can be mathematically represented by a system variable. Let ${E_1(t),\ldots,E_n(t)}$ be a set of such environmental variables, with

\begin{align*}
&\!\left\{~\begin{aligned}
\dot{E}_1 &= g_1(E_1,\ldots,E_n,V,t)\,,\quad E_1(0)=E_{1,0}\,, \\[-4pt]
\vdots & \\[-4pt]
\dot{E}_n &= g_n(E_1,\ldots,E_n,V,t) \,,\quad E_n(0)=E_{n,0}\,,
\end{aligned} \right.\, & 
\end{align*} 
being the given reaction equations. Then the ODE for the environmental stress level can be given by 
\begin{equation}
\dStress=\underbrace{\left(\sum_{j=1}^{n}\alpha^-_j\cdot \delta_j^-(E_j)\right)(1-\stress)}_{\substack{\small\mbox{increasing stress level} \\ \small \mbox{(stressful conditions)}}}\underbrace{-\left(\sum_{j=1}^{n}\alpha^+_j\cdot \delta_j^+(E_j)\right)\stress}_{\substack{\small\mbox{recovery from stress} \\ \small \mbox{(beneficial conditions)}}}\qquad\mbox{with }\stress(0)=\initStress\in[0,1]\,. \label{eq:stress}
\end{equation} 
The first term includes the dynamics raising the ESL. Since~${\stress\in[0,1]}$, the increase is bounded by one from above via the factor~${(1-\stress)}$. If the environmental factors generate beneficial growth conditions for the cells, the ESL decreases according to the second term. The parameters~$\alpha^-_j$ resp. $\alpha^+_j$ are variable-specific ``sensitivity rates'', denoting how fast the stress level increases/decreases, if the value of the associated variable~$E_j$ enters/leaves a range, which is critical for survival. With distinguishing between~$\alpha^-$ and~$\alpha^+$ per variable, it is possible to mathematically consider that cells react e.g. faster to the lack of an environmental factor, which is essential for their survival, than they recover from the deprivation once beneficial conditions are restored. Which values of~$E_j$ are considered to be critical, is given by the explicit shape of the corresponding ``influence functions''~$\delta_j^-(E_j)$ and~$\delta_j^+(E_j)$. In particular, large values of $\delta_j^-(E_j)$ promote stress, whereas large values of $\delta_j^+(E_j)$ inhibits stress.

Collecting the influences of the environment in the stress level gives rise to some practical advantages. We can choose an arbitrary number of environmental variables and include them easily into the stress reaction due to the modular structure of~\eqref{eq:stress}. With the parameters~$\alpha_j^-$ and~$\alpha_j^+$ we can consider different sensitivities of the cells to changes in different environmental factors. Although neither these parameters nor the stress level~$\stress$ can be measured directly in experiments, appropriate data of cell viability can provide enough information to estimate the values of the sensitivity rates. For this e.g. statistical methods for parameter calibration can be used, where the only prior assumption is a rough range, in which the parameter values are expected. Then, by grouping these sensitivity rates according to their magnitude, we can distinguish the corresponding environmental variables between having a fast or slow or even no impact on the cells. Under a quasi-steady state assumption, the ODE for~$\stress$ can then be simplified by choosing a time scale of interest and omitting environmental variables with no influence.

In the considered experimental setting, we have the special case of the constant nutrient saturation~$S(t)$ being the only environmental variable:
\begin{equation*}
E_1(t)=S(t)=\initS\quad \text{with}\quad g_1(E_1,V,t)=g(S,V,t)=0\,.
\end{equation*} 
We assume the ESL to increase and decrease with the same ``nutrient sensitivity rate''~$\adaptS$. The corresponding ``nutrient influence functions'' are then~$\deactFct{S}$ and $\actFct{S}$, as introduced previously in~\eqref{eq:deltas} for model~\ref{eq:ODE}. Hence, we have
\begin{equation*}
\alpha^-_1=\alpha^+_1=\adaptS\,,\quad \delta_1^-(E_1)=\deactFct{\initS}\quad \text{and}\quad \delta_1^+(E_1)=\actFct{\initS}=1-\deactFct{\initS}\,,
\end{equation*} 
i.e. given the initial condition ${\stress(0)=\initStress\in[0,1]}$ we consider the ODE
\begin{equation*}
		\dStress = \adaptS\deactFct{\initS}\cdot(1-\stress)-\adaptS\big(1-\deactFct{\initS}\big)\cdot\stress = \adaptS\big(\deactFct{\initS}-\stress\big)\,.
\end{equation*} 
Since this equation is independent of the variable~$V$, the exact solution can be calculated by separation of variables and for~$t\geq 0$ it is given as:
\begin{equation}
	{\stress(t)=\deactFct{\initS}\cdot(1-e^{-\adaptS t})+\initStress e^{-\adaptS t}}\,. \label{eq:stressSol}
\end{equation} 
Instead of the nutrient-dependent scaling functions~$\delta^\pm$ in the basic model~\ref{eq:ODE}, now the ESL influences the proliferation and starvation rate. A high stress level causes slower proliferation and faster starvation, whereas a low stress level has the opposite effect. Since by definition~${\stress(t)\in[0,1]}$, this results in the following initial value problem:
\begin{align}
	&\!\left\{~\begin{aligned}
		\dV &= \underbrace{\big(1-\stress(t)\big)\cdot\prol}_{=\,\scalProlStress(t)}\, V\left(1-\left(\frac{V}{\Kv}\right)^{m}\right)-\underbrace{\big(\nat+\stress(t)\cdot\starv\big)}_{=\,\scalDeathStress(t)} V\,, \\
		\dot\stress&=\adaptS\deactFct{\initS}\cdot(1-\stress)-\adaptS\big(1-\deactFct{\initS}\big)\cdot\stress\,, \\
		V(0)&=\initV,~ \stress(0)=\initStress\,.
	\end{aligned} \right. \tag{$\mathcal{M}_{\stress}$} & \label{eq:ODEstress}
\end{align} 
By inserting the explicit time-dependent analytical solution~\eqref{eq:stressSol}, the system reduces to a non-autonomous ODE for~$V$, which is solved numerically. We observe that model~\ref{eq:ODEstress} has an interesting relation to the model~\ref{eq:ODE}: under the assumption that critical nutrient changes immediately influence the viable cells, i.e. assuming~${\adaptS\to\infty}$, we return to model~\ref{eq:ODE}, since for~$t\geq 0$ it holds
\begin{equation*}
	\stress(t)\overset{\eqref{eq:stressSol}}{=} \deactFct{\initS}\cdot(1-e^{-\adaptS t})+\initStress e^{-\adaptS t}\overset{\adaptS\to\infty}{\longrightarrow}\deactFct{\initS}=\mbox{const.}~\Rightarrow~ \scalProlStress(t)\overset{\adaptS\to\infty}{\longrightarrow}\scaledProl~,~\scalDeathStress(t)\overset{\adaptS\to\infty}{\longrightarrow}\scaledDeath\,.
\end{equation*} 
The same result can be achieved by a quasi-steady state assumption, stating that changes in the ESL happen at a much faster time scale than changes in the tumor cell number. In this case, the stress level reaches its steady state virtually instantly, which is:
\begin{equation*}
	{\dStress=\adaptS\big(\deactFct{\initS}-\stress\big)=0\Rightarrow\stress=\deactFct{\initS}}\,.
\end{equation*}
We are aware that having multiple environmental factors included into the model would make the advantages of using the concept of the environmental stress level more evident. However, due to the lack of corresponding data, we focus on the special case of having only one environmental variable (nutrient saturation~$S$). We consider the investigations in this manuscript as a proof of principle that the environmental stress level is a feasible alternative way to model the effect of the environment on the cells.
\paragraph{Overview over all models and their mathematical properties.} The last model~\ref{eq:ODEstress} has been constructed from model~\ref{eq:ODE}, which itself uses model~\ref{eq:ODEprol} as a basis. Therefore, the number of variables and parameters increases with adding more complex dynamics. An overview over the variables and strictly positive model parameters of each system as well as abbreviating notations for important parameter terms/functions can be found in Tables~\ref{Tab:modelVars},~\ref{Tab:modelPars} and~\ref{Tab:modelNotations} below. The units of the variables and parameters are motivated by the biological setting and the measurement methods, which are explained in more detail in the following Section~\ref{subsec:experiments}.
\begin{table}[H]
	\setlength{\tabcolsep}{4pt}
	\begin{center}
		\caption{Overview over the model variables and their occurrence in the models.}
		\begin{tabular}{cccclc} \hline
			\small  \ref{eq:ODEprol} & \small \ref{eq:ODE}& \small \ref{eq:ODEstress} & \textbf{Variable} & \textbf{Meaning} & \textbf{Unit} \\\hline
			\checkmark & \checkmark & \checkmark & $V=V(t)$  & Density of viable tumor cells& \SI{}{\Cells/\milli\litre} \\
			&\checkmark & \checkmark &$S=S(t)$   & Nutrient saturation&   \mbox{\SI{}{\FBS}}\\
			&&\checkmark & $\stress=\stress(t)$  & Environmental stress level & --\\\hline
		\end{tabular}
		\label{Tab:modelVars}
	\end{center}
\end{table}
\begin{table}[H]
	\setlength{\tabcolsep}{4pt}
	\begin{center}
		\caption{Model parameters, their occurrence in the models and their meaning.}
		\begin{tabular}{cccclc} \hline
			\small  \ref{eq:ODEprol} & \small \ref{eq:ODE}& \small \ref{eq:ODEstress}  &\textbf{Parameter}  &\textbf{Meaning}& \textbf{Unit} \\\hline
			\checkmark &\checkmark &\checkmark &  $\prol$ & Maximal possible proliferation rate& \SI{}{1/\day}\\
			\checkmark &\checkmark &\checkmark &  $\Kv $ & Carrying capacity of the biological system & \SI{}{\Cells/\milli\litre}\\
			\checkmark &\checkmark &\checkmark &  $1/m$ & Strength of proliferation contact inhibition & --\\
			\checkmark &\checkmark &\checkmark & $\nat$ & Natural death rate & \SI{}{1/\day}\\
			 &\checkmark & \checkmark& $\starv$ & Maximal possible starvation rate & \SI{}{1/\day} \\
			 &\checkmark &\checkmark &  $\Sthr$& Nutrient sensitivity threshold& \mbox{\SI{}{\FBS}} \\
			 & &\checkmark & $\adaptS$ & Sensitivity rate of nutrient changes on stress level  &  \SI{}{1/\day}\\\hline
		\end{tabular}
	\label{Tab:modelPars}
	\end{center}
\end{table}
\begin{table}[H]
	\setlength{\tabcolsep}{3pt}
	\begin{center}
		\caption{Abbreviating notations for important/useful parameter terms and functions as well as their occurrence in the model equations (none of them occurs in model~\ref{eq:ODEprol}).}
		\begin{tabular}{ccrll} \hline
			\small \ref{eq:ODE}& \small \ref{eq:ODEstress}  &\textbf{Notation} & & \textbf{Meaning} \\\hline
			\checkmark &\checkmark &  $\actFct{\initS}=$ & $\frac{\initS^2}{\Sthr^2 +\initS^2}$& Influence function for nutrient-promoted dynamics\\
			\checkmark &\checkmark &   $\deactFct{\initS}=$ & $1-\actFct{\initS}$& Influence function for nutrient-inhibited dynamics \\[4pt]
			\checkmark & &   $\scaledProl=$ & $\actFct{\initS}\cdot\prol$& Nutrient-dependent net proliferation rate (constant)\\
			\checkmark & &   $\scaledDeath=$ & $\nat+\deactFct{\initS}\cdot \starv$& Nutrient-dependent net death rate (constant)\\[4pt]
			& \checkmark&  $\scalProlStress(t)=$ & $\big(1-\stress(t)\big)\cdot\prol$ & Stress-dependent net proliferation rate (varying in time)\\
			&\checkmark &  $\scalDeathStress(t)=$ & $\nat+\stress(t)\cdot \starv$ & Stress-dependent net death rate (varying in time)\\\hline
		\end{tabular}
		\label{Tab:modelNotations}
	\end{center}
\end{table}
A mathematical analysis of the models yields positivity and boundedness of the solutions, which are important features for biological reasonableness. Table~\ref{Tab:ODEfeat} summarizes the computed bounds as well as the steady states~$\equibV$ reps. ${(\equibV,\bar{\stress})^T}$ and their stability. The corresponding calculations providing these results can be found in the supplementary Section~\ref{SUPsec:mathAna}.
\begin{table}[H]
	\setlength{\tabcolsep}{4pt}
	\begin{center}
		\caption{Bounds of the solutions $V(t)$, $\stress(t)$ and steady states $\equibV$, $(\equibV,\bar{\stress})^T$ of each ODE model.}
		\begin{tabular}{c c l r l}
		\hline
		 & & \textbf{Positivity \& boundedness}  &  & \hspace*{-1.45cm}\textbf{Steady states}  \\ 
		\textbf{Model} & &  \textbf{of the solutions}  & stable & unstable \\
		\hline
		\ref{eq:ODEprol}  & & $0\leq V(t)\leq\max\left\{\initV,\,\Kv\sqrt[m]{1-\frac{\scaledDeath
			}{\scaledProl}}\,\right\}$\vphantom{$\Bigg\{$} & $\Kv\sqrt[m]{1-\frac{\scaledDeath
			}{\scaledProl}}$ & $0$ \\[10pt]
		\ref{eq:ODE} & for $\scaledDeath<\scaledProl$: & $0\leq V(t)\leq\max\left\{\initV,\,\Kv\sqrt[m]{1-\frac{\scaledDeath
			}{\scaledProl}}\,\right\}$  & $\Kv\sqrt[m]{1-\frac{\scaledDeath
		}{\scaledProl}}$&  $0$  \\
		& for $\scaledDeath\geq\scaledProl$: & $0\leq V(t)\leq\initV$  & $0$  & \\[10pt]
		\ref{eq:ODEstress} & for all $\scaledDeath,\scaledProl$: & $0\leq\,\stress(t)\leq\max\left\{\initStress,\,\deactFct{\initS}\right\}\leq 1$ & &  \\
 		 & for $\scaledDeath<\scaledProl$: & $0\leq V(t)\leq\max\left\{\initV,\,\Kv\sqrt[m]{1-\frac{\scaledDeath
 		 	}{\scaledProl}}\,\right\}$ & $\left(\begin{array}{c}
 	 	\Kv\sqrt[m]{1-\frac{\scaledDeath
 	 		}{\scaledProl}} \\ \deactFct{\initS}		\end{array}\right)$ & $\left(\begin{array}{c}
	0 \\ \deactFct{\initS}		\end{array}\right)$ \\
		& for $\scaledDeath\geq\scaledProl$: & $0\leq V(t)\leq\initV$ & $\left(\begin{array}{c}
			0 \\ \deactFct{\initS}		\end{array}\right)$ & \\
		\hline
	\end{tabular}
	\label{Tab:ODEfeat}
	\end{center}
\end{table}
\subsection{Experimental data}
\label{subsec:experiments}
We use data from the experiments described in~\cite{Lima} to calibrate the unknown parameters from Table~\ref{Tab:modelPars}. A CellTiter-Blue\textsuperscript{\textregistered} assay was used to monitor the viability of tumor cells. In particular, viable cells metabolize a provided chemical and they emit measurable light as a result of this process. Hence, the data points are fluorescence intensity measurements. This way, viability was measured once every day and there are four biological replicates of each measurement to estimate statistical significance and repeatability.

Let $I^\text{total}$ be an intensity measurement of a specific cell line at time point~$\tilde{t}$. Excluding the corresponding background intensity~$\iBG$ of the cell-free medium, the fluorescence intensity produced by viable cells~$I^V$ is assumed to be directly proportional to the density of viable tumor cells~${V(\,\tilde{t}\,)}$. The experiments were performed with different initial cell densities between~$\num{e4}$ and $\num{e5}$ cells per milliliter, which motivates the unit of~$V$ and eventually leads to the relation
\begin{equation*}
	I^V=I^\text{total}-\iBG\propto V(\,\tilde{t}\,) ~\mbox{\SI{}{\Cells\per\milli\litre}}\qquad\Rightarrow\qquad n=\frac{I^V}{V(\,\tilde{t}\,)}\in(0,1)\,,
\end{equation*} 
where $n$~denotes the proportionality constant translating fluorescence intensity to cell density. Whenever we refer to ``intensity measurements'' in the following, we mean the fluorescence produced by the cells~$I^V$ and neglect the superscript~``$V$'' for better readability.

\begin{figure}[H]
	\begin{center}
		\includegraphics[width=.85\textwidth]{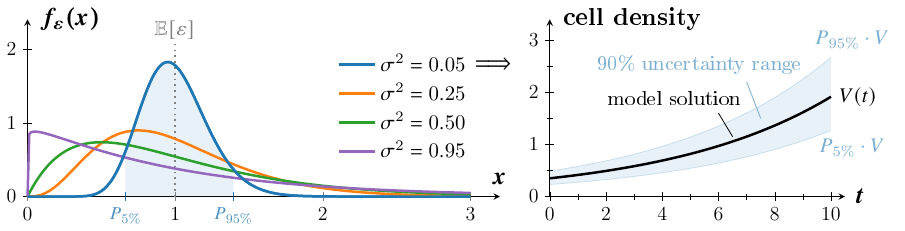}
		\caption{Left: Probability density function~$f_\varepsilon$ of ${\varepsilon\sim\Gamma\left(1/\sigma^2\,,1/\sigma^2\right)}$ for varying variances~$\sigma^2$ and the percentiles~$P_{5\%}$ and~$P_{95\%}$ for ${\sigma^2=0.05}$. Right: Resulting 90\% uncertainty range (shaded area) around an exemplary model solution~$V(t)$}
		\label{Fig:Gamma}
	\end{center}
\end{figure}

For nutrition the cells are supplemented with a particular concentration between $0\%$ to $10\%$ of fetal bovine serum (FBS). A supplementation with \SI{0}{\fbs} does not provide the cells with any nutrients, whereas \SI{10}{\fbs} generates optimal growth conditions. The nutrient supply is kept constant throughout the whole duration of each experiment, i.e. the nutrient saturation is assumed to stay at its initial level ${S_0\in\left[0,1\right]}$ at any time. Overall, this motivates:
\begin{equation*}
	S(t)=S(0)=S_0~\forall\,t\geq 0\quad\mbox{with unit}\quad[S]=\mbox{\SI{}{\FBS}}.
\end{equation*} 
The remaining variable, the environmental stress level~$\stress=\stress(t)$, is an auxiliary variable and especially immeasurable. Hence, it has no experimental counterpart and is used as a dimensionless quantity.
\paragraph{How are uncertainties considered?} In reality, the equation~${n=I/V(\,\tilde{t}\,)}$ is not rigorously fulfilled. This can be due to e.g. model inadequacy or biological fluctuations of the cells' metabolism, which affect measurement accuracy. To capture this uncertainty, we assume a multiplicative noise for each element of a set of~${M\in\mathbb{N}}$ measurements~${\{I_i\}_{i=1}^{M}}$ and a set of model solutions~${\{V_i\}_{i=1}^{M}}$ considering the corresponding values of~$t$, $\initV$, and $\initS$. In the Bayesian framework that we adopt (see Section \ref{subsec:calibration}) we model this to be a random variable~$\varepsilon_i$ such that for~${i=1,\ldots,M}$:
\begin{equation}
	I_i=n\, V_i\cdot\varepsilon_i\,.  \label{eq:intensity} 
\end{equation} 
Let~$\varepsilon_i$ be i.i.d. and have the unimodal and continuous distribution of a random variable~$\varepsilon$ with probability density function~(PDF)~$f_\varepsilon$. Then, the following properties should hold:
	\begin{align*}
	[\mbox{P1}]:~~&\text{supp}(f_\varepsilon)\subseteq\IR^+\,,\mbox{i.e. all measurements are positive}; \\
	[\mbox{P2}]:~~&\mean[\varepsilon]= 1\,,\mbox{i.e. measurements are accurate on average};\\[-4pt]
	[\mbox{P3}]:~~&f_\varepsilon(x)\overset{x\to \infty}{\longrightarrow}0 \mbox{ and } f_\varepsilon(x)\overset{x\to 0}{\longrightarrow}0\,,\mbox{i.e. outliers are possible but not likely}.
	\end{align*} 
Different distributions are possible to accomplish these properties. A small number of shape parameters and an easy calculation to ensure property~[P2] are desirable. Therefore, we choose a Gamma distribution~${\varepsilon\sim\Gamma(a,b)}$, ${a,b>0}$ with a few restrictions. This is a plausible choice for multiplicative noise and often used in imaging theory (see e.g.~\cite{gamma1,gamma2,gamma3}). Property~[P1] and the desired behavior near infinity of~[P3] are satisfied by definition. The corresponding PDF is given by ${f^{a,b}_\varepsilon(x)=\frac{b^a}{\Gamma(a)}\,x^{a-1}e^{-bx}}$, where ${\Gamma(\,\cdot\,)}$ denotes the Gamma function. To fulfill the remaining properties, we observe
\begin{equation*}
	\mean[\varepsilon]=\frac{a}{b}\overset{\mbox{\tiny [P2]}}{=}1\,\mbox{, if }~b=a\qquad \Rightarrow~\lim\limits_{x\to 0}f^{a,a}_\varepsilon(x)=\lim\limits_{x\to 0}\left(\frac{a^a}{\Gamma(a)}\,x^{a-1}e^{-ax}\right)\overset{\mbox{\tiny [P3]}}{=}0\,\mbox{, if }~a>1\,.
\end{equation*} 
Note that by constraining~${b=a>1}$ the shape of the distribution depends only on the parameter~$a$. In fact, $a$ is directly related to the standard deviation~$\sigma$ and hence the variance of the distribution: ${\sigma^2=\Var(\varepsilon)=\frac{a}{a^2}=\frac{1}{a}\Rightarrow a=\frac{1}{\sigma^2}}$.
Therefore, for~${i=1,\ldots,M}$ the uncertainty factor~$\varepsilon_i$ for a particular measurement~$I_i$ can be modeled by
\begin{equation}
\varepsilon_i=\frac{I_i}{n\cdot V_i}\sim\Gamma\left(\frac{1}{\sigma^2},\frac{1}{\sigma^2}\right)\mbox{ with }\sigma^2\in(0,1)\,. \label{eq:noise}
\end{equation} 
We use the percentiles~$P_{5\%}$ and~$P_{95\%}$ of the Gamma distributed uncertainty factors~$\varepsilon_i$ to define the
\begin{equation}
	\text{``90\% (uncertainty) range'' around the solution }V:\quad \big[V\cdot P_{5\%}\,,\, V\cdot P_{95\%}\big]\,.\label{eq:range}
\end{equation} 
In particular, given a specific noise variance~$\sigma^2$, the model expects 90\% of the measurements within this interval, whereas respectively 5\% are expected below and above it. The left side of Figure~\ref{Fig:Gamma} depicts exemplary plots of the PDF~$f_\varepsilon$ for different values of~$\sigma^2$. It also shows the positive skewness of the distribution. Under consideration of the measurement method, this is a reasonable feature for the uncertainty factors, assuming the cells might not metabolize the assay to their full potential. On the right side of Figure~\ref{Fig:Gamma} we see an example of a 90\% uncertainty range around a solution.
The reason for considering a multiplicative (and, for instance, not additive) noise term is twofold. From a mathematical perspective, it allows to preserve positivity of the data. A more practical motivation is the reasonable assumption that the fluorescence noise of the intensity measurements is proportional to the density of viable cells. This is also supported by the observation of a larger variance for experiments with larger cell numbers in our data. Furthermore, multiplicative noise has demonstrated to be better suited for fluorescence than additive noise in other experimental settings~\cite{Vsavsik02}.
\paragraph{How can the measurements be utilized?} One of the experiments in~\cite{Lima} monitors the dependence between cell viability and nutrient supply by providing five different nutrient concentrations~$\initS$ over a period of 7 days: this results in five data sets {``D1''--``D5''}. These measurements are employed to calibrate the unknown model parameters for models~\ref{eq:ODE} and~\ref{eq:ODEstress} with Bayesian inversion methods (for details see following Section~\ref{subsec:calibration}). Another experiment in~\cite{Lima} was designed to investigate cell behavior under optimal growth conditions, i.e. the cells were provided with~\SI{10}{\fbs} over a period of 21 days. We denote the data set generated by this experiment with ``D6''. It is used additionally to {D1--D5} to validate the calibration results with model~\ref{eq:ODEprol}. All experiments start with several populations of different initial size~$\initV$. Table~\ref{Tab:experiments} shows the experimental setup of all data sets and the resulting initial conditions in the models.

\begin{table}[H]
	\setlength{\tabcolsep}{4pt}
	\begin{center}
		\caption{Overview over the data sets and the corresponding initial values. D1--D5 are used to calibrate the models~\ref{eq:ODE} and~\ref{eq:ODEstress} and D6 for validation of the results with model~\ref{eq:ODEprol}.}
		\begin{tabular}{ccclcc} \hline
			&\textbf{Nutrition}& \textbf{Duration}& \makebox[0pt][l]{\textbf{Initial values in the models}} & & \\[-2pt]
			\textbf{Data set} &(in \SI{}{\fbs})& (in days) &$\boldsymbol{\initV}$ & $\boldsymbol{\initS}$ &$\boldsymbol{\initStress}$ \\\hline
			\textbf{D1} & 10.0 & 7 & $1.00, 0.50, 0.25$ & $1.00$& $0.00$  \\
			\textbf{D2}& \hphantom{0}7.5 & 7  &$1.00, 0.50, 0.25$ & $0.75$ & $0.00$\\
			\textbf{D3}& \hphantom{0}5.0 & 7  & $1.00, 0.50, 0.25$ & $0.50$ & $0.00$  \\
			\textbf{D4}& \hphantom{0}2.5 & 7 &$1.00, 0.50, 0.25$ & $0.25$ & $0.00$\\
			\textbf{D5}& \hphantom{0}0.0 &7  &$1.00, 0.50, 0.25$ & $0.00$ & $0.00$  \\[4pt]
			\textbf{D6} & 10.0 & 21& $1.00, 0.50, 0.25, 0.10, 0.05$ & ($1.00$)& ($0.00$) \\\hline
		\end{tabular}
		\label{Tab:experiments}
	\end{center}
\end{table}
\subsection{Parameter calibration: Bayesian inversion and Sequential Monte Carlo (SMC)}
\label{subsec:calibration}
The task to identify the unknown true parameters from given data is called the ``inverse problem''. We solve it using a Bayesian approach, which leads, under mild assumptions, to a naturally well-posed inverse problem~\cite{Stuart10}. Furthermore, this approach allows to quantify the uncertainty in the estimated parameters and hence it brings more information than a deterministic inversion method, at the price of a higher computational cost.

We collect all parameters to be estimated in a vector ${\theta\in\Theta\subseteq \mathbb{R}^d}$ (${d\in\mathbb{N}}$), where $\Theta$ denotes the parameter space. Assuming the considered data consist of $M$ intensity measurements (excluding background intensity), we collect them in a vector ${\mathcal{I}=\big({I}_1,\,\ldots\,,\,{I}_M\big)^T=\big({I}_i\big)^M_{i=1}\in\mathbb{R}^M}$. Defining ${\mathcal{G}: \Theta\rightarrow \mathbb{R}^M}$ with ${\theta\mapsto\big(\mathcal{G}_i(\theta)\big)_{i=1}^M}$ as the forward operator mapping parameter values to the corresponding intensities, such a measurement can be rewritten as
\begin{equation}\label{eq:forward}
	{I}_i \overset{\eqref{eq:intensity}}{=} n\,V_i\cdot\varepsilon_i=\mathcal{G}_i(\theta)\cdot\varepsilon_i\,,
\end{equation} 
where $V_i$ is the corresponding model solution to \eqref{eq:ODEprol}, \eqref{eq:ODE} or \eqref{eq:ODEstress} using $\theta$ as parameters and $\varepsilon_i$ is the multiplicative noise. Note that we use parts of the data sets D1--D5 for parameter estimation, hence each measurement~${I}_i$ can refer to a different nutrient condition, initial cell density, and time point. Therefore, $V_i$ has to be calculated in consideration of the corresponding values for~$\initS$, $\initV$ and $t$.

Now, we want to consider all measurements in~${\mathcal{I}=\big({I}_i\big)^M_{i=1}}$ collectively. The parameter vector~$\theta$ and the noise~${\big(\varepsilon_i\big)^M_{i=1}}$ are modeled as multi-dimensional random variables taking values in $\Theta$ and $\mathbb{R}^M$ respectively. The Bayesian formulation of the problem is the following: Given a prior (measure) $\mu_0$ on $\Theta$, compute the posterior (measure) $\mu^{\mathcal{I}}$ given the data~$\mathcal{I}$. The prior in the Bayesian setting is the correspondent of a regularization in deterministic inverse problems~\cite{Stuart10} and it reflects the knowledge about the parameters before including any information given by the data, whereas the posterior describes the knowledge after seeing the data. Let $\pi_0$ and $\pi^{\mathcal{I}}$ denote the probability densities of $\mu_0$ and $\mu^{\mathcal{I}}$, respectively. By Bayes' formula, we have
\begin{equation}\label{eq:posterior}
	\pi^{\mathcal{I}}(\theta) = \frac{L(\mathcal{I}\,|\,\theta)\cdot\pi_0(\theta)}{\int_{\Theta}L(\mathcal{I}\,|\,\theta)\cdot\pi_0(\theta)\,\mbox{d}\theta}\propto L(\mathcal{I}\,|\,\theta)\,\pi_0(\theta)\,,
\end{equation} 
where $L$ is the data likelihood~\cite{Stuart10}. The proportionality constant of relation~\eqref{eq:posterior} depends only on $\mathcal{I}$. It is called ``model evidence'' and it can be used to quantitatively compare two models (see Section~\ref{subsec:modelComp}). We remind from the previous Section~\ref{subsec:experiments} that~\eqref{eq:noise} states i.i.d. ${\varepsilon_i\sim\Gamma\left(\frac{1}{\sigma^2}, \frac{1}{\sigma^2}\right)}$ with ${\sigma^2\in(0,1)}$ for every measurement~${I}_i$,~${i=1,\ldots,M}$. Using this together with~\eqref{eq:forward}, the data likelihood of ${\mathcal{I}=\big({I}_i\big)^M_{i=1}}$ is:
\begin{equation*}
	L(\mathcal{I}\,|\,\theta)\propto\prod_{i=1}^M \left(\frac{{I}_i}{\mathcal{G}_i(\theta)}\right)^{1/\sigma^2-1}\exp\left(-\frac{1}{\sigma^2}\frac{{I}_i}{\mathcal{G}_i(\theta)}\right)\,.
\end{equation*} 
\paragraph{How to sample from the posterior?} In order to make predictions, we want to sample from the given posterior distribution \eqref{eq:posterior}. However, this has a complicated, concentrated density, so we cannot sample from it exactly with a random number generator. To approximate the posterior measure, we use therefore the Sequential Monte Carlo (SMC) method, which we explain now based on~\cite{KBJ14}. In SMC one considers a sequence of intermediate distributions ${(\mu_k)_{k=0}^N}$\,, such that $\mu_0$ is the prior and ${\mu_N=\mu^{\mathcal{I}}}$ coincides with the posterior distribution. The probability density~$\pi_k$ of the intermediate measure $\mu_k$ can be defined by
\begin{equation}\label{eq:pit}
	\pi_k(\theta) = \frac{1}{Z_k} \prod_{i=1}^k L_i({I}_i\,|\,\theta)\, \pi_0(\theta)\quad\mbox{or, equivalently,}\quad \pi_k(\theta) =\frac{1}{\tilde{Z}_k} L_k(I_k\,|\,\theta)\, \pi_{k-1}(\theta)\,,
\end{equation} 
where $Z_k$ and $\tilde{Z}_k$ are normalizing constants and
\begin{equation*}
	L_i({I}_i\,|\,\theta)  = \frac{1}{Z_i}\left(\frac{{I}_i}{\mathcal{G}_i(\theta)}\right)^{1/\sigma^2-1}\exp\left(-\frac{1}{\sigma^2}\frac{{I}_i}{\mathcal{G}_i(\theta)}\right)
\end{equation*} 
is the likelihood associated to observation~${I}_i$ with normalization constant~$Z_i$,~${i=1,\ldots,M}$. The intermediate densities could also be constructed with an adaptive approach using tempering~\cite{propKernel}. However, since the considered data measures the quantity of interest in a time series, the presented filtering method is computationally more efficient. The SMC algorithm samples sequentially from the intermediate measures $\mu_k$ using a weighted swarm of samples, called particles. Let ${P\in\mathbb{N}}$ be the sample size, i.e. the number of particles. At the $k$-th iteration (${k=1,\ldots,N}$) the algorithm leads to a collection of particles ${\big\{\theta_p\big\}_{p=1}^{P}}$ with associated weights ${\big\{W^k_p\big\}_{p=1}^{P}}$, which gives the approximation
\begin{equation*}
	\pi_k(\theta)\approx \sum_{p=1}^{P} W^k_p\delta_{\theta_p}(\theta)\,.
\end{equation*} 
The SMC algorithm is summarized in Algorithm \ref{alg:SMC}: we achieve appropriately weighted particles to approximate~${\mu_N=\mu^{\mathcal{I}}}$ by starting with uniformly weighted particles distributed according to the prior~$\mu_0$ (line~\ref{algline:start}) and iteratively move the samples from the previous measure~$\mu_{k-1}$ to $\mu_k$ in a selection~(lines~\ref{algline:sel1}-\ref{algline:sel2}) and a mutation step~(line~\ref{algline:MCMC})~\cite[Ch. 5]{D04}, which are explained in detail in the following paragraphs.
\begin{algorithm}[H]
	\caption{Sequential Monte Carlo}\label{alg:SMC}
	\begin{algorithmic}[1]
		\State $k=0$: sample $\theta_p\sim\mu_0$ and set $W^0_p=\frac{1}{P}\,,~~p=1,\ldots,P$\label{algline:start}
		\For{$k=1,\ldots,N$}
		\State $w_p^k =  L_k(I_k\,|\,\theta)\, W^{k-1}_p\,,~~  W_p^k = w_p^k/\big(\sum_{p=1}^{P} w_p^k\big)\,,~~ p=1,\ldots,P$ \Comment{selection step}\label{algline:sel1}
		\If{$P_\text{eff}<P_{\text{thr}}$} resample:
		\State $(i)$: sample indices acc. to distribution~$\mathcal{R}$ of particle indices: $\big(\pi_1,\ldots,\pi_{P}\big)\sim\mathcal{R}\big(W_1^k,\ldots,W_{P}^k\big)$
		\State $(ii)$: set $\theta_p = \theta_{\pi_p}$ and $W^k_p=\frac{1}{P}\,,~~p=1,\ldots,P$
		\EndIf \label{algline:sel2}
		\State move $\theta_p\sim\kappa_k(\theta_p\,,\,\cdot\,)\,,~~p=1,\ldots,P$\Comment{$\mu_k$-invariant mutation step}\label{algline:MCMC}
		\EndFor
	\end{algorithmic}
\end{algorithm}

\textsl{Selection step}. We start with a collection of particles ${\big\{\theta_p\big\}_{p=1}^{P}}$, distributed according to $\mu_{k-1}$. Their weights~${\big\{W^{k-1}_p\big\}_{p=1}^{P}}$ are updated to ${\big\{W^{k}_p\big\}_{p=1}^{P}}$ by importance sampling: for ${p=1,\ldots,P}$ we have
\begin{equation*}
	W_p^k = \frac{w_p^k}{\sum_{p=1}^{P} w_p^k}\quad\mbox{with}\quad
	w_p^k \overset{\eqref{eq:pit}}{=}  L_k(I_k\,|\,\theta)\,W^{k-1}_p\,.
\end{equation*}
We see that $W_p^k$ is normalized to ensure that ${\sum_{p=1}^{P}W_p^k=1}$, i.e. having a probability distribution. Note that the importance sampling only changes the weights and not the particles. However, if there are many particles with low weights, the estimation is only as accurate as a Monte Carlo approximation with a very small number of particles~\cite{pendu}. In this case, the reweighing step is followed by a resampling step, where the particles are replaced according to their updated weights. Resampling is needed, if the effective sample size
\begin{equation*}
	P_\text{eff}\approx \left(\sum_{p=1}^{P} \left(W_p^k\right)^2\right)^{-1}
\end{equation*} 
is small, which we check by comparison with a threshold ${P_{\text{thr}}=\tau P\,,\,\tau\in(0,1)}$:
\begin{equation*}
	P_\text{eff}\,\begin{cases}
		\,< P_{\text{thr}} & \Rightarrow \mbox{ resample } \big\{\theta_p\big\}_{p=1}^{P} \mbox{ according to }\big\{W^{k}_p\big\}_{p=1}^{P} \mbox{ and uniformly weigh new particles}\,,\\
		\,\geq P_{\text{thr}} & \Rightarrow \mbox{ do not resample}\,. \\
	\end{cases}
\end{equation*} 
This discards particles with low weight and improves the representation of the distribution~$\mu_k$.

\smallskip

\textsl{Mutation step}. Performing only selection steps will eventually lead to degeneracy in the diversity of the particle population. In particular, after some resampling steps, few particles will survive and be replicated. Therefore, we introduce diversity in the particles by moving them according to a Markov Chain Monte Carlo (MCMC) kernel ${\kappa_k(\,\cdot\,,\,\cdot\,)}$. This kernel is $\mu_k$-invariant, i.e. it does not modify the particle distribution. We adopt the adaptive strategy developed in~\cite{propKernel} to construct such a MCMC kernel. A random walk Metropolis-Hastings (MH) proposal is used on each univariate component, conditionally independently. More precisely, remembering that each particle is a vector ${\theta_p\in\mathbb{R}^d}$, MH proposes ${\big\{q_p\big\}_{p=1}^{P}}$ with ${q_p\in\mathbb{R}^d}$ by computing
\begin{equation*}
	(q_p)_j = (\theta_p)_j + \epsilon^k_j\cdot\xi_j\quad\mbox{with}\quad \xi_j\sim\mathcal{N}(0,1)\,,
\end{equation*}
where ${(\,\cdot\,)_j}$ denotes the $j$-th component (${j=1,\ldots,d}$) and the scale~$\epsilon^k_j$ is tuned to the acceptance rate of the previous SMC iteration. This means, we choose $\epsilon^k_j$ adaptively as
\begin{equation*}
	\epsilon^k_j = \rho_k \sqrt{\widehat{\text{Var}}\Big((\theta)_j\Big)}~\mbox{with recursively defined scaling parameter}~ \rho_k=\begin{cases}
		\rho_{k-1}\cdot 2 &\text{if }a_{k-1}>0.30,\\
		\rho_{k-1}\,/\,2&\text{if }a_{k-1}<0.15,\\
		\rho_{k-1}& \text{otherwise},
	\end{cases}
\end{equation*} 
where $a_{k-1}$ is the average acceptance rate over the particles at the previous iteration and ${\widehat{\text{Var}}\Big((\theta)_j\Big)}$ denotes the empirical marginal variance from the $j$-th components of all particles. It is also possible to adapt the proposal using empirical covariances instead~\cite{adaptMH}, which would be computationally more expensive in view of possible extensions to parameters, which are random fields and therefore very high dimensional. Eventually, the particles are moved by randomly accepting the proposed ones:
\begin{equation*}
	\mbox{set }\theta_p = q_p\quad\text{ with probability }\quad \min\left\{\frac{\pi_k(q_p)}{\pi_k(\theta_p)}\,,\,1\right\}.
\end{equation*} 
To improve the mixing, it is possible to repeat this process more than once by applying~$\kappa_k$ again on the moved particles. Afterwards, the total ratio of accepted particles determines the $a_k$, which is used to scale~$\epsilon^{k+1}_j$ in the MCMC update of the next iteration.
\paragraph{How is the theory applied to the models?} Parameter estimation is performed in Python adapting the code provided in~\cite{pendu}. The set of parameters, which need to be calibrated, can be distinguished between model parameters (see Table~\ref{Tab:modelPars}) and hyper parameters. Latter includes the unknown uncertainty variance~$\sigma^2$ and the proportionality constant~$n$ from relation~\eqref{eq:noise}. All calibrations are performed using the SMC method with sample size~${P=50\,000}$ using the resampling threshold~${P_\text{thr}=75\%\,P}$. In the mutation step, five MCMC updates are performed, i.e. line~\ref{algline:MCMC} in Algorithm~\ref{alg:SMC} is performed five times.

\begin{figure}[H]
	\begin{center}
		\includegraphics[width=0.95\textwidth]{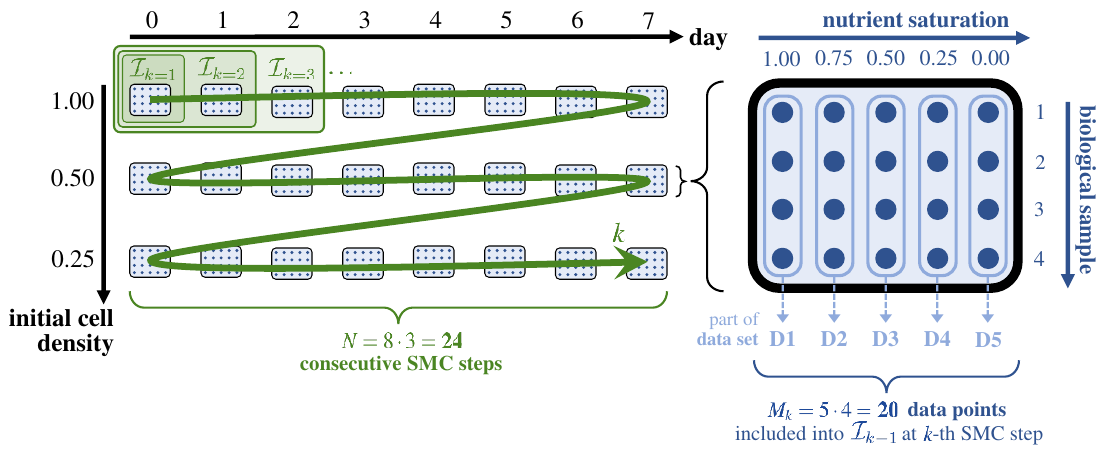}
		\caption{Schematic description of how the data sets~D1--D5 are used to calibrate the parameters step by step with SMC}
		\label{Fig:MCsteps}
	\end{center}
\end{figure}

The parameters are calibrated over the course of ${N=8\cdot 3=24}$~SMC~steps: As depicted on the left side of Figure~\ref{Fig:MCsteps}, starting with the data corresponding to ${\initV=1.00}$ and ${t=0}$, the SMC steps iterate over the eight points in time (inner loop) and over the three separate seeding densities (outer loop). With each step~$k$, another set of ${5 \cdot 4 =20}$~data points, containing all measurements of D1--D5 at a specific time point~$t$ regarding a particular initial cell density~$\initV$ (see right side of Figure~\ref{Fig:MCsteps}), is included into the considered data collection~$\mathcal{I}_k$ for the selection step (line~\ref{algline:sel1} in Algorithm~\ref{alg:SMC}). This maximizes the information about the effect of varying nutrients on the cells at each calibration step. In particular, in the first step a set of measurements~${\mathcal{I}_1=\big\{\mathcal{I}_{1,i}\}^{M_1}_{i=1}}$ with ${M_1=5 \cdot 4 =20}$ is considered and with each SMC step ${20}$ more data points are included incrementally: ${\mathcal{I}_1\subset\mathcal{I}_2\subset\ldots\subset\mathcal{I}_N}$ with ${|\mathcal{I}_k|=M_k=k\cdot 20}$, ${k=1,\ldots,N}$. In total, the parameters are calibrated using ${M_N=24\cdot 20=480}$ measurements.
We want to construct informative prior distributions for the parameters. Information about the parameters' magnitude can be described by using an appropriate uniform distribution~${\unif(a,b)}$ on an interval~${[a,b]}$. If we can additionally assume that the neighborhood of a parameter value~$h$ has a high probability, a triangular distribution~${\Triang(a,h,b)}$ on~${[a,b]}$ is used, where ${h\in[a,b]}$ is the mode of the distribution, i.e. the value that is most likely to be sampled. Note that the mode can also lie on the interval bounds, if low/large values within the interval are assumed to have high probability. We set ${h=(a+b)/2}$ to the center of the interval, if we expect the borders ${a,b}$ to have a low probability but there is no particular tendency to a value within the interval.

The doubling time~$T$ of the used cell type is larger than one day while in exponential growth~\cite{prol1, prol2}, i.e. ${T>1}$. We can estimate an upper bound for the growth rate~$\prol$:
\begin{equation*}
2\initV=V(T)=\initV e^{\prol T}\quad\Rightarrow~ 2= e^{\prol T} \quad\Rightarrow~ \prol=\frac{\ln(2)}{T}\overset{T>1}{<}\ln(2)<1\,.
\end{equation*} 
Therefore, it is reasonable to expect ${\prol\in(0,1)}$. Under optimal nutrient conditions, the cell population size increases. This translates to the parameter relation ${\prol>\nat}$, which can also be written as ${\nat=c_1\cdot\prol}$ with ${c_1\in(0,1)}$. Since the parameter~$c_1$ is bounded, we calibrate this one instead of~$\nat$. A similar reparametrization can be done for the starvation rate~$\starv$\,: it is reasonable to assume that the cells die faster from starvation than from natural causes in case of a nutrient-free environment, which leads to ${\starv>\nat\Rightarrow\starv=\frac{\nat}{c_2}=\frac{c_1}{c_2}\prol}$ with ${c_2\in(0,1)}$. Furthermore, the underlying experimental setting motivates the assumption that the population size does not exceed ${3\times\num{e5}}$ cells per milliliter. All cell lines are seeded in initial densities below the carrying capacity and ${\initV=1}$ is the largest seeding density, i.e. ${\Kv\in(1,3)}$ is plausible. The restrictions ${m>1}$ and ${\Sthr\in(0,1})$ with ${\Sthr\ll 1}$ are motivated by the modeling framework. Utilizing all this information, we adopt the following prior distributions for the model parameters of~\eqref{eq:ODE}:
\begin{alignat*}{3}
	\prol&\sim\unif(0,1)\,,\qquad c_1&&\sim\Triang(0,1/2,1)\,,\qquad c_2&&\sim\Triang(0,1/2,1)\,, \\
	\Kv&\sim\unif(1,3)\,,\qquad m&&\sim\unif(1,12)\,,\hspace*{0.6cm} \text{and}\hspace*{0.65cm} \Sthr&&\sim\Triang(0,0,1)\,.
\end{alignat*}
The same prior distributions are used for model~\ref{eq:ODEstress}. Its equations only have one additional model parameter~$\adaptS$, for which we do not have any particular information -- we set its prior distribution to
\begin{equation*}
\adaptS\sim\unif(0,12)\,.
\end{equation*} 
To consider a certain degree of confidence in the measurements, we assume a small uncertainty variance~${\sigma^2\in(0,1)}$ is more likely than a large one. Therefore, we use ${\Triang(0,0,1/2)}$ as its prior. All experiments are started with cells from a batch with optimal nutrient conditions. If they are put into a nutrient-free environment without going through a weaning process beforehand, they can undergo a starvation shock. This might disturb or decrease the cells' ability to metabolize the chemical for the fluorescence measurements. To consider this in the hyper parameters, we allow the data sets D1--D4 (${\initS>0}$) and D5 (${\initS=0}$) to have different uncertainty variances~${\sigma^2_{\text{D1:4}}}$ resp.~${\sigma^2_{\text{D5}}}$ and proportionality constants~$n_{\text{D1:4}}$ resp.~$n_{\text{D5}}$, where potentially ${n_{\text{D1:4}}\geq n_{\text{D5}}}$. Using the reparametrization~${n_{\text{D5}}=c_n\cdot n_{\text{D1:4}}}$ with ${c_n\in(0,1)}$, leads to the prior distributions
\begin{equation*}
	\sigma^2_{\text{D1:4}}\,,\, \sigma^2_{\text{D5}}\sim\Triang(0,0,1/2)\,,\quad 
	n_{\text{D1:4}}\sim\unif(0,1/2)\,,\quad\mbox{and}\quad c_n\sim\Triang(0,1,1)\,.
\end{equation*} 
A large value for the uncertainty variance~$\sigma^2$ allows a larger deviation of the model solution from the data. For the purpose of model comparison, using the same uncertainty variance for both models is desirable to increase comparability. Hence, we first calibrate each model separately to get an estimate for their variances~$\sigma^2_{\text{D1:4}}$ and~$\sigma^2_{\text{D5}}$. Then, we take the average~$\bar{\sigma}^2$ of the means respectively over both models, i.e.
\begin{equation}
	\bar \sigma^2_{\text{D1:4}}=\frac{1}{2}\left(\, \mean_S\left[\sigma^2_{\text{D1:4}}\right]+\mean_\stress\left[\sigma^2_{\text{D1:4}}\right]\,\right)\qquad\mbox{and}\qquad\bar \sigma^2_{\text{D5}}=\frac{1}{2}\left(\,\mean_S\left[\sigma^2_{\text{D5}}\right]+\mean_\stress\left[\sigma^2_{\text{D5}}\right]\,\right)\,, \label{eq:avgNoise}
\end{equation} 
where ``$S$'' resp. ``$\stress$'' in the subscripts of the expected value~$\mean$ indicate the underlying model~\ref{eq:ODE} resp.~\ref{eq:ODEstress}, which is used for the calibration. These average values~$\bar \sigma^2_{\text{D1:4}}$ and $\bar \sigma^2_{\text{D5}}$ are then used deterministically and are especially not estimated anymore in further calibrations. For better comparability, we start each SMC algorithm from the same prior particle sample. In particular, the algorithm is performed for model~\ref{eq:ODEstress} and the generated initial sample but without the component regarding to parameter~$\adaptS$, is used to start the calibration of model~\ref{eq:ODE}.
\subsection{Model comparison}
\label{subsec:modelComp}
\begin{figure}[H]
	\begin{center}
		\includegraphics[width=0.6\textwidth]{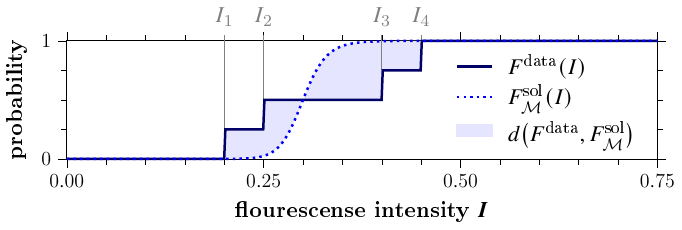}
		\caption{Plots of the data distribution function~$F^{\text{data}}$ given by four exemplary intensity measurements~${I_1,\ldots,I_4}$, an approximate prediction distribution function~$F^{\text{sol}}_\mathcal{M}$ and the corresponding validation metric~${d\big(F^{\text{data}},F^{\text{sol}}_\mathcal{M}\big)}$ given by the enclosed area between the two functions}
		\label{Fig:valid}
	\end{center}
\end{figure}

The calibration results of each model are compared in a quantitative manner. For this, we calculate different comparison measures and validate the model solutions with the data sets~D1--D6. An overview over the applied methods is given in the following paragraphs of this section.
\paragraph{How can the quality of fit be quantified?} We use the validation metric proposed in~\cite{validMetric} to compare the model prediction at a given point in time with the corresponding set of measurements. Their mismatch is measured as the area between the data distribution~$F^{\text{data}}$ and the prediction distribution~$F^{\text{sol}}_\mathcal{M}$ using the calibration results of model~$\mathcal{M}$, mathematically defined by the metric
\begin{equation}
d\big(F^{\text{data}},F^{\text{sol}}_\mathcal{M}\big)=\int_{0}^{\infty}\left|\,F^{\text{data}}(I)-F^{\text{sol}}_\mathcal{M}(I)\,\right|~\mbox{d}I\,.\label{eq:valid}
\end{equation} 
For a set~${\{\,I_{i}\,\}_{i=1,\,\ldots\,,\,M}}$ of $M$~intensity measurements, the data distribution function is given by
\begin{equation*}
F^{\text{data}}(I)=\sum_{i=1}^M\mathbb{I}\Big( I_i,\,I\Big)\quad\mbox{with}\quad \mathbb{I}\Big( I_i,\,I\Big)=\begin{cases}
1/M &\mbox{for } I_i\leq I\,, \\
0 &\mbox{for } I_i> I\,.
\end{cases}
\end{equation*}
The prediction distribution function~$F^{\text{sol}}_\mathcal{M}$ is determined approximately: for each weighted particle~$\theta_p$ of the posterior, we calculate the empirical cumulative distribution function from the solution of model~$\mathcal{M}$. We get a set~${\{\,\tilde{I}_{p}\,\}_{p=1,\,\ldots\,,\,P}}$ with~${\tilde I_{p}=n_p\cdot V_p}$, where $V_p$ is the model solution using the parameter sample~$\theta_p$ scaled with the corresponding sample~$n_p$ of the proportionality constant. The prediction distribution is then 
\begin{equation*}
F^{\text{sol}}_\mathcal{M}(I)\approx\sum_{p=1}^P\mathbb{I}_W\Big(\tilde I_{p},\,I\Big)\quad\mbox{with}\quad \mathbb{I}_W\Big(\tilde I_{p},\,I\Big)=\begin{cases}
W_{p}^N &\mbox{for }\tilde I_{p}\leq I\,, \\
0 &\mbox{for }\tilde I_{p}> I\,,
\end{cases}
\end{equation*} 
where~$W^N_p$ is the final weight of the $p$-th particle after the SMC algorithm is finished. Figure~\ref{Fig:valid} depicts an exemplary plot of the above distribution functions and the resulting validation metric for a set of four measurements.
\paragraph{How can the different model approaches be compared?}
We use the Bayes factor~\cite{Bayes} for a quantitative model comparison. It gives the ratio between the respective model evidences. Latter describe the posterior probability of the data given the model type in consideration of the parameter distribution. For a model~$\mathcal{M}$, its evidence is defined by the marginal likelihood of a set of measurements~$\mathcal{I}$:
\begin{equation*}
p(\mathcal{I}\,|\,\mathcal{M})=\int_{\Theta}\pi_0(\theta\,|\,\mathcal{M})\cdot L(\mathcal{I}\,|\,\theta,\mathcal{M})\,\mbox{d}\theta\,,
\end{equation*} 
where ${\pi_0(\theta\,|\,\mathcal{M})}$ denotes the prior on $\theta$ using model~$\mathcal{M}$ and ${L(\mathcal{I}\,|\,\mathcal{M},\theta)}$ the likelihood of~$\mathcal{I}$ given~$\mathcal{M}$ and its parameters~$\theta$. Following~\cite{evid}, the evidence~$Z_k$ at the $k$-th SMC step can be approximated by
\begin{equation*}
Z_k\approx Z_0\cdot \prod_{j=1}^{k}\sum_{p=1}^{P}\frac{L_j(\mathcal{I}_{j}\,|\,\theta_p,\mathcal{M})}{L_{j-1}(\mathcal{I}_{j-1}\,|\,\theta_p,\mathcal{M})}\cdot W_p^{j-1}\,,
\end{equation*} 
where $P$ is the number of particles, ${L_j(\mathcal{I}_{j}\,|\,\theta_p,\mathcal{M})}$ is the likelihood of observing all data sets included until calibration step~$j$ given the parameter sample~$\theta_p$ with model~$\mathcal{M}$, and $ W_p^{j-1}$ denotes the normalized weight of~$\theta_p$ at step~${j-1}$. Since we start the calibrations with priors, from which we can sample directly, the initial evidence is~${Z_0=1}$. The Bayes factor of two models~$\mathcal{M}_1$ and $\mathcal{M}_2$ regarding the same set of measurements~$\mathcal{I}$ is then given as
\begin{equation*}
\mathcal{Z}_{1:2}=\frac{p(\mathcal{I}\,|\,\mathcal{M}_1)}{p(\mathcal{I}\,|\,\mathcal{M}_2)}\,.
\end{equation*} 
For  ${p(\mathcal{I}\,|\,\mathcal{M}_1)>p(\mathcal{I}\,|\,\mathcal{M}_2)}$, i.e. ${\mathcal{Z}_{1:2}>1=10^0}$\,, the strength of evidence can be described by the following scale~\cite{evid}:
\begin{equation}
	\log_{10}\big(\mathcal{Z}_{1:2}\big)\in\,\begin{cases}
		~\big(0,\frac{1}{2}\big] &\text{barely worth mentioning}\,, \\
		~\big(\frac{1}{2},1\big] &\text{substantial support for }\mathcal{M}_1\,, \\
		~\big(1,2\big] &\text{strong support for }\mathcal{M}_1\,, \\
		~\big(2,\infty\big) &\text{decisive support for }\mathcal{M}_1\,. \\
	\end{cases}\label{eq:scaleZ}
\end{equation} 
\section{Results and Discussion}
\label{sec:results}
This chapter investigates the results of the model calibrations. We compare the corresponding estimated parameters of both models~\ref{eq:ODE} and~\ref{eq:ODEstress} and quantify the quality of fit to the data (Section~\ref{subsec:resComp}). The calibration results can be used to identify parameter correlations (Section~\ref{subsec:resCorr}) and to study the reaction of the cells to different nutrient concentrations (Section~\ref{subsec:resNut}). Additionally to the calibration data sets, the models are validated with further data (Section~\ref{subsec:resVal}).

To investigate dispersion of the results, we perform 12 runs of the SMC algorithm. In the following, the numerical deviations are given in terms of the mean and the 95\% prediction/confidence interval of a normal distribution~${\mathcal{N}(\mu,\sigma^2)}$\,, i.e. denoted by ${\mu\pm 1.96\sigma}$. If not indicated otherwise, this formula is also used for error bars in visual representations of the deviations.

According to the equations in~\eqref{eq:avgNoise}, a precalibration of the models yields the uncertainty variances
\begin{equation*}
\bar \sigma^2_{\text{D1:4}}=0.0355\qquad\mbox{and}\qquad\bar \sigma^2_{\text{D5}}=0.2410\,,
\end{equation*} 
which are used deterministically for further calibrations. We see that for the data sets~D1--D4~(${\initS>0}$) the variance is considerably smaller than for data set~D5~(${\initS=0}$). On the one hand, this could be a consequence of the previously described starvation shock, which might disturb the cells' metabolism. On the other hand, data set D5 contains mainly measurements of low cell density, for which the measurement accuracy might be decreased because of weak fluorescence.
\subsection{Comparison of the model calibration results}
\label{subsec:resComp}
Regarding the estimated posteriors, the ESS is larger than $P_\text{thr}$ after performing the last SMC step for all model calibrations. This implies that the final approximated posterior distribution is nearly as accurate as sampling directly from the correct probability measure. We calculate the expected values and variances (see Tables~\ref{Tab:deviations1} and~\ref{Tab:deviations2}) of the model parameters' marginal posteriors.
\begin{table}[H]
	\setlength{\tabcolsep}{4pt}
	\begin{center}
		\caption{Expected values~$\mean(\,\cdot\,)$ of the posterior distributions for each model parameter.}
		\begin{tabular}{c c c c c}
			\hline
			~& ~ & ~ & ~ & \llap{\parbox{8.65cm}{\textbf{Calibrated model parameters}}} \\[-2pt]
			\textbf{Model} & $\boldsymbol{\prol}$ & $\boldsymbol{\nat}$ & $\boldsymbol{\starv}$ & $\boldsymbol{\adaptS}$ \\
			\hline
			\ref{eq:ODE} & $0.435\pm0.043$ & $0.103\pm0.042$ & $0.186\pm0.031$ & ~ \\
			\ref{eq:ODEstress} & $0.437\pm0.023$ & $0.106\pm0.025$ & $0.196\pm0.018$ & $6.930\pm2.561$ \\[4pt]
			~ & $\boldsymbol{\Kv}$ & $\boldsymbol{m}$ & $\boldsymbol{\Sthr}$ &\\
			\cline{2-5}
			\ref{eq:ODE} & $1.740\pm0.133$ & $4.731\pm2.371$ & $0.104\pm0.011$ &\\
			\ref{eq:ODEstress} & 
			$1.731\pm0.098$ & $5.315\pm2.964$ & $0.106\pm0.007$ &\\
			\hline
		\end{tabular}
		\label{Tab:deviations1}
	\end{center}
\end{table}
\vspace*{-5pt}
\begin{table}[H]
	\setlength{\tabcolsep}{4pt}
	\begin{center}
		\caption{Variances~$\Var(\,\cdot\,)$ of the posterior distributions for each model parameter.}
		\begin{tabular}{c c c c c}
			\hline
			~& ~ & ~ & ~ & \llap{\parbox{8.65cm}{\textbf{Calibrated model parameters}}} \\[-2pt]
			\textbf{Model} & $\boldsymbol{\prol}$ & $\boldsymbol{\nat}$ & $\boldsymbol{\starv}$ & $\boldsymbol{\adaptS}$ \\
			\hline
			\ref{eq:ODE} & 
			$0.007\pm0.017$ & $0.006\pm0.014$ & $0.010\pm0.020$ & ~  \\
			\ref{eq:ODEstress} & $0.010\pm0.022$ & $0.009\pm0.020$ & $0.012\pm0.024$ & $3.211\pm3.014$ \\[4pt]
			~ & $\boldsymbol{\Kv}$ & $\boldsymbol{m}$ & $\boldsymbol{\Sthr}$ &\\
			\cline{2-5}
			\ref{eq:ODE} & 
			$0.028\pm0.061$ & $0.751\pm0.777$ & $0.004\pm0.009$  &\\
			\ref{eq:ODEstress}& $0.033\pm0.067$ & $1.071\pm1.782$ & $0.004\pm0.008$ &\\
			\hline
		\end{tabular}
		\label{Tab:deviations2}
	\end{center}
\end{table}
The tables show very similar values for both models and the differences between them lie within the range of the numerical variations of the SMC algorithm. The resulting average parameter values can be concluded as biologically reasonable, since we chose a prior distribution ensuring this. Up to some adaptions regarding units and reparametrization, the reasonableness of the parameters is also supported by comparing their values to the calibration results of~\cite{Lima}, where a set of similar models is used together with parts of the same data.

Taking a closer look at the nutrient sensitivity rate~$\adaptS$, we observe a large average variance of its marginal posterior (${\Var(\adaptS)\approx3.211}$) as well as large numerical deviations of its expected value (${\pm\,2.561}$) and variance (${\pm\,3.014}$). A possible reason for these large values is that only a minority of the measurements contain meaningful information to estimate~$\adaptS$. On the one hand, this parameter influences the cell growth mainly for early measurements, when the cells react to the nutrient change from the batch colony (${S=1}$) to the provided nutrient supply of the respective experiment (${S=\initS\in[0,1]}$). On the other hand, the experiments with ${\initS\gg\Sthr}$ show little measurable reaction of the cells to nutrient changes due to the small value of the nutrient threshold~$\Sthr$, and therefore lack information to estimate~$\adaptS$. Overall, approximately less than 10\% of the measurements contain meaningful information to estimate~$\adaptS$, which leads to a higher uncertainty in its estimation. Nevertheless, we can argue in Section~\ref{subsec:resNut} that this does not have a big influence on the solution of model~\ref{eq:ODEstress}.

\begin{figure}[H]
	\begin{center}
		\includegraphics[width=.9\textwidth]{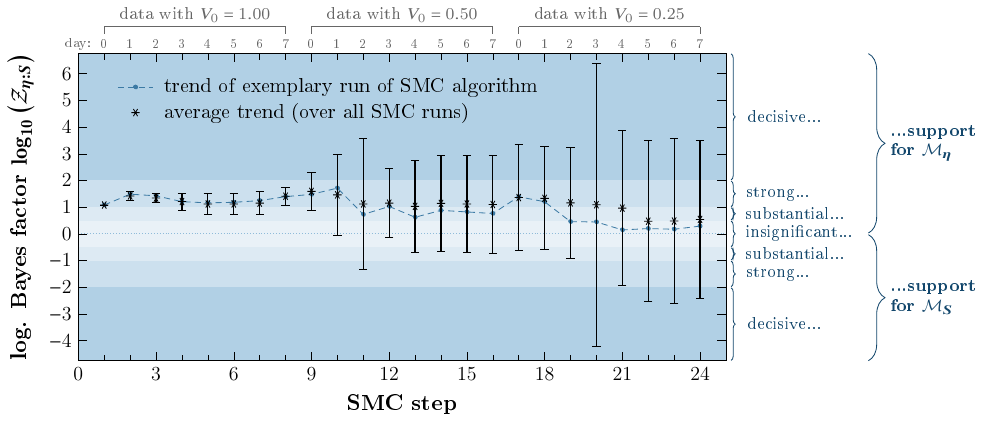}
		\caption{Logarithm of the Bayes factor~$\mathcal{Z}_{\stress:S}$ over the course of the SMC steps, where the error bars indicate the 95\% confidence interval resulting from the various runs of the SMC algorithm. The dashed trendline shows the evolution for a single run, which estimated the model parameters close to the average values of all runs. The gray labels on top of the plot give the included data points at each SMC step}
		\label{Fig:evid}
	\end{center}
\end{figure}

We investigate the ratio of the resulting values of the validation metric~\eqref{eq:valid} of each model. This does not show any preference of a particular model (see Table~\ref{Tab:validResSUP} in supplementary Section~\ref{SUPsec:calRes}).
To get a quantitative statement for the models' quality of fit to the data, we calculate the logarithm of the Bayes factor: ${\log_{10}\big(\mathcal{Z}_{\stress:S}\big)=\log_{10}\left(\frac{p(\mathcal{I}\,|\,\mathcal{M}_\stress)}{p(\mathcal{I}\,|\,\mathcal{M}_S)}\right)}$. In Figure~\ref{Fig:evid} the average trend of this value over the course of model calibration steps is shown and interpreted according to the scale given in \eqref{eq:scaleZ}. Additionally to the averaged trend, the plot shows the evolution of the Bayes factor for one particular run of the algorithm (dashed line). For this run, the resulting expected values of the parameters' posteriors match approximately with the ones in Table~\ref{Tab:deviations1}. We will refer to the corresponding posteriors later in Section~\ref{subsec:resCorr}, when we investigate the correlation between the model parameters.
For the first eight calibration steps, the Bayes factor indicates strong evidence to prefer model~\ref{eq:ODEstress} over~\ref{eq:ODE}. Until this point only data regarding~${\initV=1.00}$ are considered. Incremental inclusion of the next data set (${\initV=0.50}$) weakens the weight of the support for model~\ref{eq:ODEstress}, but still allows to conclude a tendency of the evidence towards this model. This can be observed further until SMC~step~19. At step~20 (i.e. with inclusion of data with ${t\approx 3}$ and ${\initV=0.25}$), the uncertainty suddenly increases drastically, which does not allow a clear interpretation of the Bayes factor anymore. The 95\% confidence interval reaches from areas of decisive support for model~\ref{eq:ODEstress} to decisive support for model~\ref{eq:ODE}. After this step the uncertainty decreases again but stays on a high level, still not allowing a clear interpretation of the Bayes factor. Overall, the average trend of the Bayes factor indicates support for model~\ref{eq:ODEstress}, despite needing an additional variable and parameter compared to model~\ref{eq:ODE}, which is already implicitly penalized by the Bayes factor~\cite{Bayes}. The modeling advantages of using the environmental stress level presented in Section~\ref{subsec:model} still hold true, even though the present experimental data are not best-suited for showing quantitatively in a prominent way the benefits of this approach. Further investigations with better-suited data are needed to demonstrate the full potential of the environmental stress level.

\begin{figure}[H]
	\begin{center}
		\includegraphics[width=.7\textwidth]{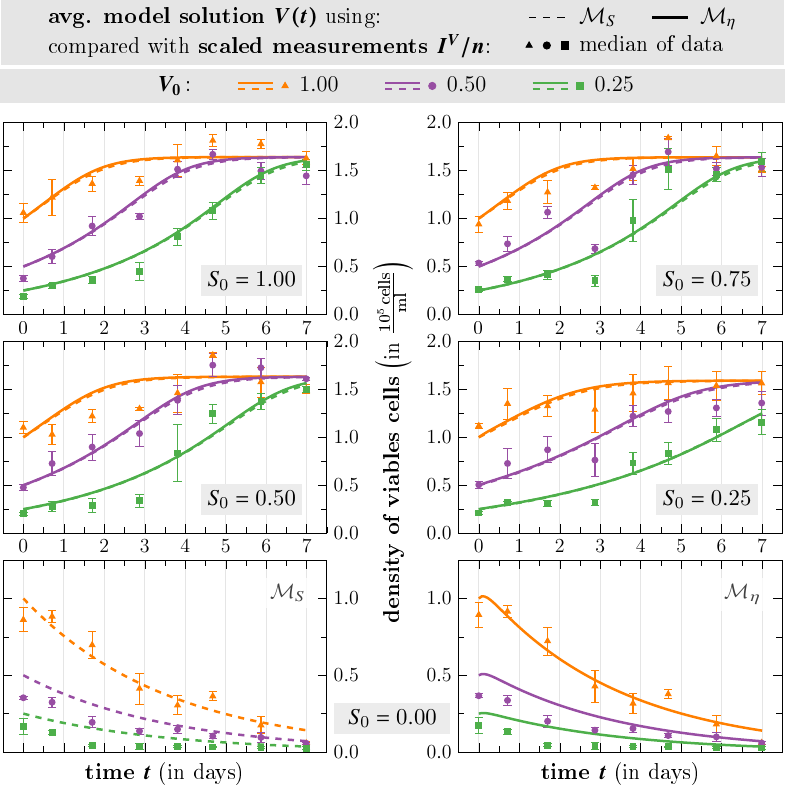}
		\caption{Time evolution of the average model solution~$V$ of model~\ref{eq:ODE} resp.~\ref{eq:ODEstress} for different~$\initV$ compared to the corresponding median of the scaled measurements~${I^V/n}$. The error bars show the median absolute deviation considering the four biological replicates}
		\label{Fig:resFull}
	\end{center}
\end{figure}

To understand the large deviations in the Bayes factor, especially the one at SMC step 20, we take a closer look at the model solution and the data. We use the average expected values from Table~\ref{Tab:deviations1} (model parameters) and supplementary Table~\ref{Tab:deviationsSUP} (scaling factors $n_\text{D1:4}$ and $n_\text{D5}$) to calculate the corresponding average model solution~$V(t)$ and scale the measurements~$I^V$. The resulting time evolution of~$V(t)$ in comparison with the calibration data sets~D1--D5 can be seen in Figure~\ref{Fig:resFull}.
In general, we observe a good fit for both models. As expected, the similar model parameter values in Table~\ref{Tab:deviations1} result in nearly the same solution for each model -- at least for all cases of~${\initS>0}$ (first two rows of plots). In these cases the estimated proportionality constant is ${n_{\text{D1:4}}\approx 0.24}$ for both models, whereas for~${\initS=0}$ they differ: ${n_{\text{D5}}\approx 0.19}$ resp. ${n_{\text{D5}}\approx 0.18}$ (\ref{eq:ODE} resp. \ref{eq:ODEstress}). This results in different scaling of data set~D5 (bottom row of Figure~\ref{Fig:resFull}). However, the difference of the value of~$n_\text{D5}$ between the models lies in the scale of numerical deviations of the SMC algorithm (see supplemented Table~\ref{Tab:deviationsSUP}).

Taking a closer look into the measurements in Figure~\ref{Fig:resFull}, the measured population size drops for all~$\initS$ at the third day of the experiment (${t\approx 3}$). This measurement might be an outlier and it corresponds to the data included in SMC~step~20, where we observe a large deviation in the Bayes factor. To investigate if this drop influences the quality of the fit, we consider the underlying measurement uncertainty. Dependent on the uncertainty variance~$\sigma^2$, we can calculate the percentiles~$P_{5\%}$ and~$P_{95\%}$ of the Gamma distributed uncertainty factor~$\varepsilon$:
\begin{align*}
	\bar \sigma^2_{\text{D1:4}}=0.0355\quad&\Rightarrow\quad P_{5\%}=0.712\,,~P_{95\%}=1.329\,,\\
	\bar \sigma^2_{\text{D5}}=0.2410\quad&\Rightarrow\quad P_{5\%}=0.350\,,~P_{95\%}=1.920\,.
\end{align*}
These determine the 90\% uncertainty range~\eqref{eq:range} for the calibration data D1--D4 resp. D5, i.e. give an interval around the solution~$V(t)$, where the model expects 90\% of the measurements.

We want to know how many measurements are actually situated where they are expected. For each data set, we count the scaled measurements~${I^V/n}$, which are situated below/within/above the 90\% uncertainty range of the corresponding average solution~$V$. A graphic overview over the calculated percentages can be found in the supplementary Section~\ref{SUPsec:calRes} (Figure~\ref{Fig:resNoise90}). We see that averaging over the data sets, for both models roughly 90\% of the scaled data is actually situated within the 90\% uncertainty range. Regarding the remaining data points, we observe that the model solutions tend to be larger than the measurements, since about 8.3--8.5\% of the data points lie below the 90\%~range, whereas only about 2.1\% are above. To investigate this observation, we take a closer look at the data below the 90\% range by checking if the discrepancy focuses on specific measurements, see Table~\ref{Tab:noise5}.
\begin{table}[H]
	\setlength{\tabcolsep}{4pt}
	\begin{center}
		\caption{Percentage of data points in D1--D5, which are below the 90\% uncertainty range for different initial cell densities~$\initV$ (rows) and days~$t$ (columns). If there are numerical deviations between the runs, they are given in terms of the 95\% confidence interval of a normal distribution.}
		\begin{tabular}{cccccccccc} \hline
			~ & & & & & & & & & \llap{\parbox{12.25cm}{\textbf{Day of measurement (i.e. $\boldsymbol t$)}}} \\[-2pt]
			\textbf{Model~~} &\hspace*{-0.4cm} $\boldsymbol{\initV}$ & \textbf{0} & \textbf{1} & \textbf{2} & \textbf{3} & \textbf{4} & \textbf{5} & \textbf{6} & \textbf{7}\\\hline
			\ref{eq:ODE} & \hspace*{-0.3cm}$\boldsymbol{1.00}$ &	$0	$&	$0.8\pm 	3.8$&	$0$&	$0$&	$0$&	$0$&	$0$&	$0$\\
			&\hspace*{-0.4cm} $\boldsymbol{0.50}$& \hphantom{0}$5.8\pm	3.8$&$0$&	\hphantom{0}$5.8\pm	3.8$&	$35.0$&	$0$&$0$&$0$&	$0$ \\
				&\hspace*{-0.4cm} $\boldsymbol{0.25}$ & $13.8\pm	6.1$&	$5.0$&	$30.8\pm	3.8$&	$73.8\pm	4.4$&	$5.4\pm	2.8$&	$0$&	$1.3\pm	4.4$&	$4.2\pm	3.8$\\[5pt]
			\ref{eq:ODEstress} & $\hspace*{-0.3cm}\boldsymbol{1.00}$ &
			$0$&	$0$&	\hphantom{0}$0.8\pm	3.8$&	\hphantom{0}$5.0$&	$3.3\pm	4.8$&	$0$&	$5.0$&	$5.0$\\
			&\hspace*{-0.4cm} $\boldsymbol{0.50}$&	$5.0$&	$0$&	\hphantom{0}$7.9\pm	5.0$&	$35.0$&	$0$&	$0$&	$0$&	$0$ \\
			&\hspace*{-0.4cm} $\boldsymbol{0.25}$ &	$14.2\pm	3.8$&$5.0$& $32.1\pm	5.0$&	$72.9\pm	5.0$&	$5.4\pm	2.8$&	$0$&	$1.3\pm	4.4$&	$5.0$\\
			\hline
		\end{tabular}
		\label{Tab:noise5}
	\end{center}
\end{table}
On the third day, the measurements show an extraordinarily large amount of data below the 90\% uncertainty range, where only 5\% would be expected. This supports the hypothesis that the measurements on that day might be outliers. The percentages are higher the smaller~$\initV$ is: approx. 35\% resp. 68--78\% of the measurements are below the 90\% range for ${\initV=0.50}$ resp. ${\initV=0.25}$. This negative correlation is expectable, since the width of the range~${\big[V\cdot P_{5\%}\,,\, V\cdot P_{95\%}\big]}$ decreases with smaller~$V$ (see right side of Figure~\ref{Fig:Gamma}). This also explains another observation: within the measurements of a particular day, more data points tend to be below the 90\% range for smaller~$\initV$. Overall, these discrepancies between measurements and expectation can be a reason for the high deviations of the Bayes factor when including data with ${\initV\leq 0.50}$ in the SMC calibration (see Figure~\ref{Fig:evid}). In particular, the sudden increase in deviation at SMC step 20 can be linked to the potential outliers for ${t\approx 3}$ and ${\initV=0.25}$.

The decrease of population size on the third day can be observed for almost all seeding densities~$\initV$ and nutrient saturations~$\initS>0$, although there is no biological reason for this as long as there are nutrients available. Hence, it is reasonable to assume that the discrepancies result from an unknown experimental bias rather than modeling inaccuracy or biological variation.
\subsection{Correlations between the model parameters}
\label{subsec:resCorr}
\begin{figure}[H]
	\begin{center}
		\includegraphics[width=.8\textwidth]{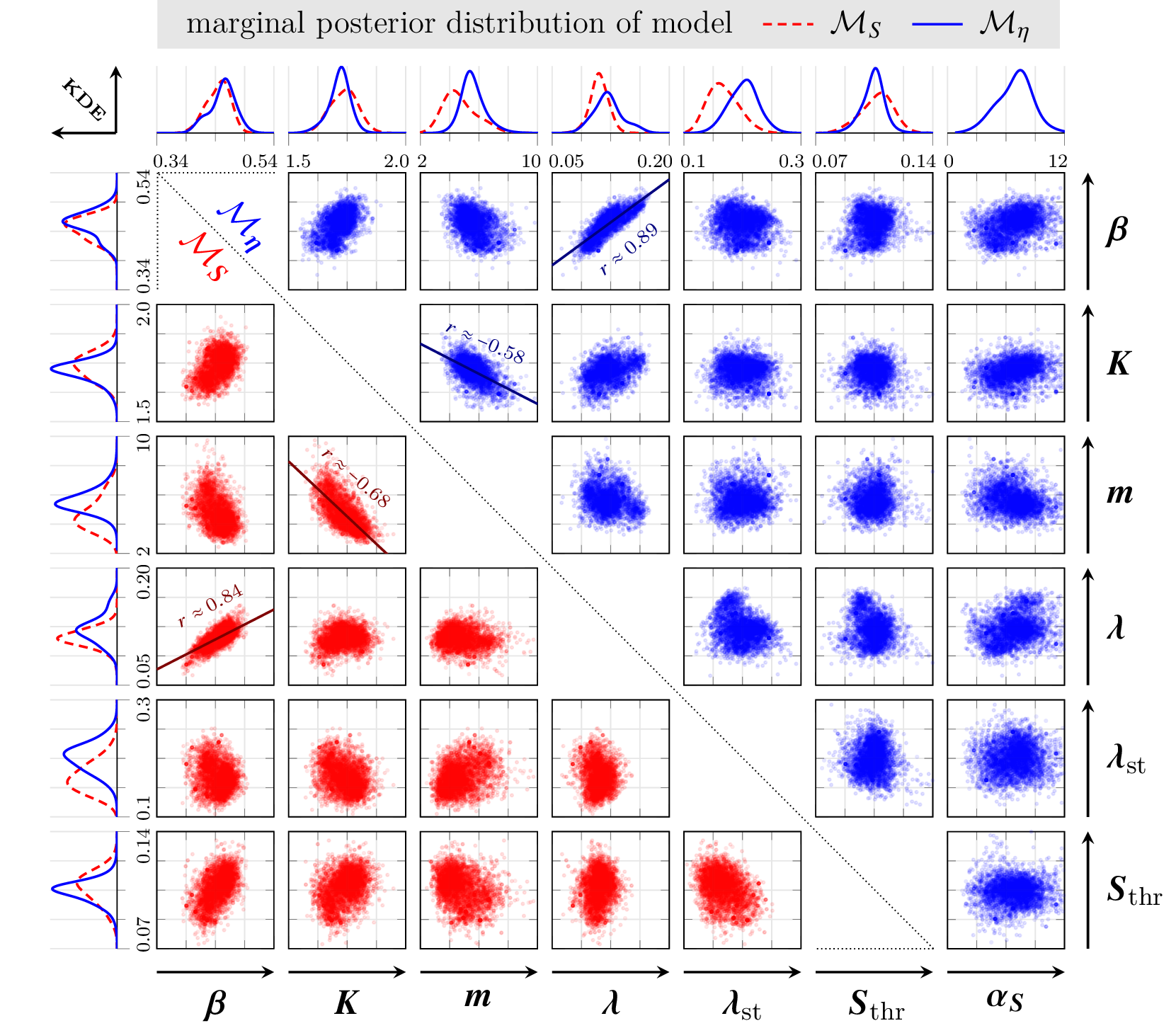}
		\caption{The horizontal/vertical axis labels are given on the bottom/right side. Plots in first row/column: Marginal posterior distributions of the model parameters using model~\ref{eq:ODE} resp.~\ref{eq:ODEstress} (note: magnitude of KDE axis is neglected and not the same for the different parameters). Remaining plots: Scatter plots for $5000$~samples drawn from the 2D distributions of corresponding pairwise parameter combinations (red, below diagonal: \ref{eq:ODE}; blue, above diagonal: \ref{eq:ODEstress}). A regression line and correlation coefficient~$r$ is depicted if a linear correlation is observable}
		\label{Fig:corrMat}
	\end{center}
\end{figure}

To identify parameter correlations, we investigate the estimated posteriors from a single run of the SMC algorithm. We choose a run whose resulting expected values of the parameters' posteriors match approximately the ones in Table~\ref{Tab:deviations1}. Figure~\ref{Fig:corrMat} shows the corresponding marginal posteriors of the model parameters for each model (first row/column) and samples drawn from the 2D posterior distributions of pairwise parameter combinations. The shape of the scatter plots can give an impression regarding the correlation of parameters.
It is observable that the SMC algorithm leads to unimodal distributions. The scatter plots indicate no correlation for most of the parameters, as the points are not concentrated along a curve. In particular, the stress-related parameter~$\adaptS$ does not correlate with other parameters (last column of scatter plots). For both models, a linear correlation can be observed between~$\prol$ and~$\nat$ (\ref{eq:ODE}: ${r\approx 0.84}$; \ref{eq:ODEstress}: ${r\approx 0.89}$) as well as between~$\Kv$ and~$m$ (\ref{eq:ODE}: ${r\approx -0.68}$; \ref{eq:ODEstress}: ${r\approx -0.58}$). These correlations are statistically significant, since their p-values are sufficiently small (${\ll0.01}$).

We have seen in Section~\ref{subsec:model} that the ratio between the growth rate~$\prol$ and the natural death rate~$\nat$ is crucial for the behavior of the population: it determines the ``net'' growth rate ${\netProl=\prol-\nat}$ and the ``net'' capacity ${\netKv=K(1-\nat/\prol)^{\frac{1}{m}}}$. Since ${\netProl,\,\netKv>0}$ should hold, those quantities require the relation~${\prol>\nat}$. This motivates the observed strong positive correlation (${r>0.7}$). For the carrying capacity~$\Kv$ and the parameter~$m$, we see a moderate negative correlation (${-0.7<r<-0.5}$). These parameters regulate the proliferation contact inhibition in the ODEs of both models~\ref{eq:ODE} and~\ref{eq:ODEstress} with the term~$(V/\Kv)^m$. This term is smaller, the larger~$\Kv$ or~$m$ are. Hence, these parameters need to be negatively correlated to mathematically describe a certain level of contact inhibition.
\subsection{Influence of different nutrient concentrations on the cell behavior}
\label{subsec:resNut}
To affect cell growth and death, both models use the same influence functions~$\delta^\pm(S)$. In model~\ref{eq:ODE} these functions scale the rates~$\prol$ and~$\starv$ directly, whereas in model~\ref{eq:ODEstress} they scale them indirectly via the stress level~$\stress(t)$. Therefore, the value of the nutrient sensitivity threshold~$\Sthr$ should not depend on the choice of the model. This behavior is actually observable in the results: The expected value of the estimated nutrient sensitivity threshold~$\Sthr$ appears to be nearly identical for both models with ${\mean\left[\Sthr\right]\approx 0.1}$ (see Table~\ref{Tab:deviations1}), resulting in the same shape of the influence function~$\deactFct{S}$.

\begin{figure}[H]
	\begin{center}
		\includegraphics[width=0.75\textwidth]{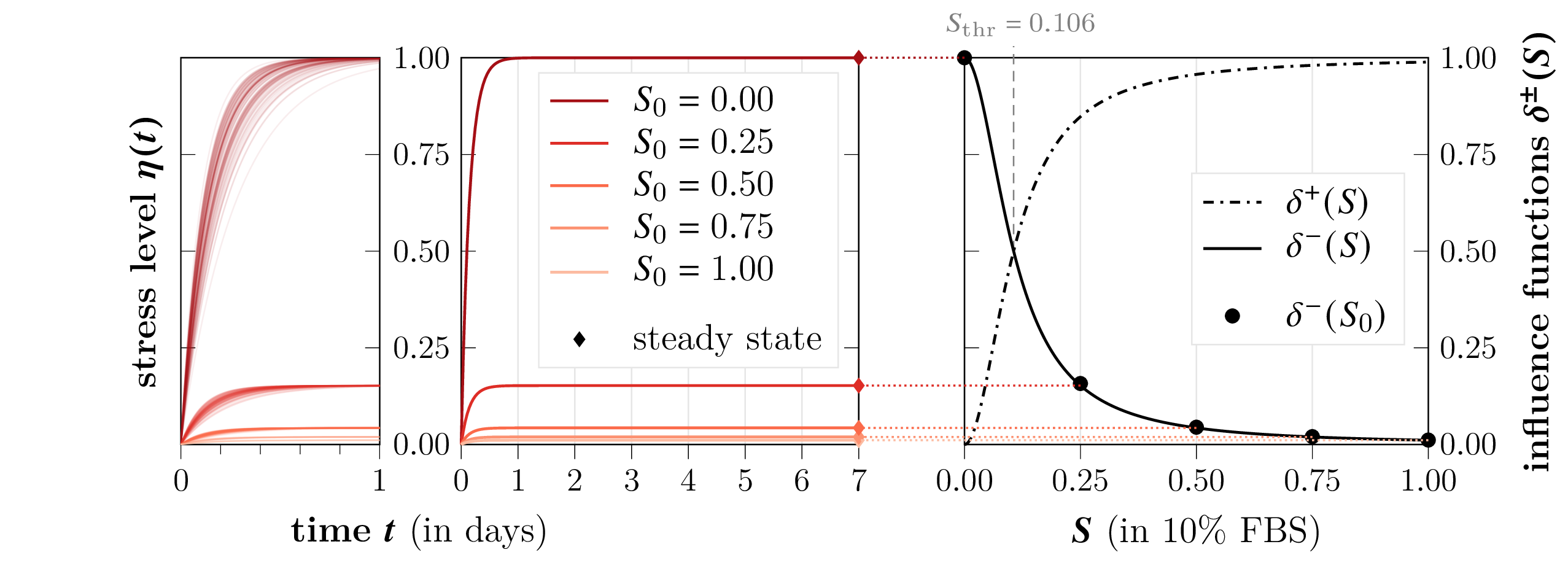}
		\caption{Time evolution of the environmental stress level~$\stress(t)$ with ${\Sthr=0.106}$ for different nutrient conditions~$\initS$ (left plot: ${\adaptS\sim\mathcal{N}(6.930,2.567/1.96)}$; middle plot: ${\adaptS=6.930}$), where its steady states coincide with the corresponding values~$\deactFct{\initS}$ of the influence function (right plot)}
		\label{Fig:stress}
	\end{center}
\end{figure}

We use the corresponding expected values from Table~\ref{Tab:deviations1} to calculate the stress level in model~\ref{eq:ODEstress}. The resulting time evolution of the stress level~$\stress(t)$ is depicted in the first two plots of Figure~\ref{Fig:stress}. The left plot shows~$\stress(t)$ for fixed $\Sthr$ and varying~$\adaptS$, while the middle one also fixes~$\adaptS$ to its calibrated mean. Fixing~$\Sthr$ is feasible, as we concluded in previous Section~\ref{subsec:resCorr} that it is uncorrelated to~$\adaptS$. The third plot shows the corresponding influence functions~$\delta^\pm(S)$.
Since the estimated value of the sensitivity parameter~$\adaptS$ is large, the steady state is already reached after approximately one day. This holds despite the relatively large deviations of the calibrated value for~$\adaptS$ (see left plot). Therefore, variation of~$\adaptS$ within the calibrated scale does not change the solution of model~\ref{eq:ODEstress} significantly. As calculated in Section~\ref{subsec:model}, the steady state of~$\stress(t)$ coincides with the function value~$\delta^-(\initS)$ of the influence function (see right plot). Hence, the high nutrient sensitivity explains why the calibrated model parameters and the resulting model solutions turn out to be very similar.
\subsection{Model validation with additional data and ``limit model''~\ref{eq:ODEprol}}
\label{subsec:resVal}
\begin{figure}[H]
	\begin{center}
		\includegraphics[width=.88\textwidth]{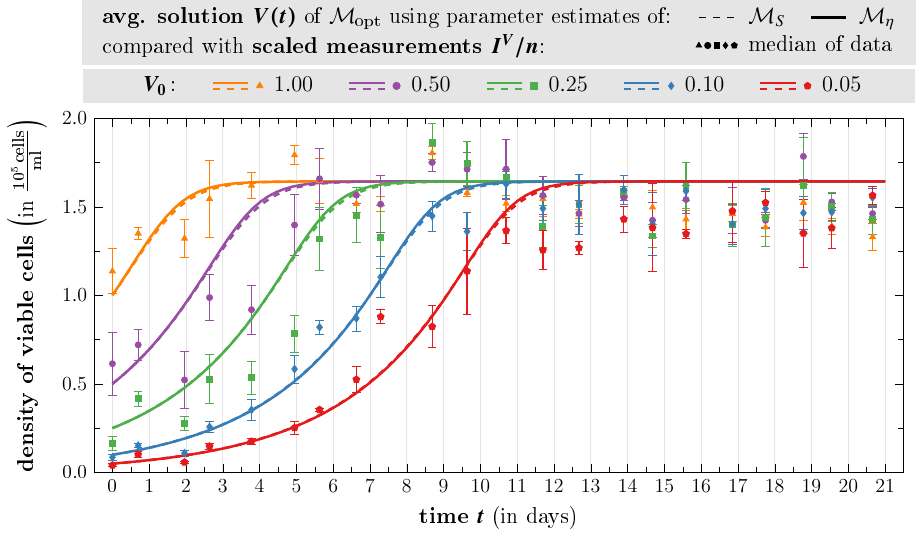}
		\caption{Time evolution of the average model solution~$V$ of model~\ref{eq:ODEprol} for different~$\initV$ using the parameter estimates from the calibration of model~\ref{eq:ODE} resp.~\ref{eq:ODEstress} compared to the corresponding median of the scaled measurements~${I^V/0.24}$. The error bars show the median absolute deviation considering the four biological replicates}
		\label{Fig:resOpt}
	\end{center}
\end{figure}

For additional validation, we calculate the average solution of the ``limit model''~\ref{eq:ODEprol} using the posteriors (determined with data sets~D1--D5) to compare it with data set~D6. In this experiment tumor cells were seeded in five different initial densities and supplied with 10\%~FBS for the duration of 21~days (see Table~\ref{Tab:experiments}). Hence, it holds~${\initS=1>0}$ (like in data set~D1) and we use ${\mean[n_{\text{D1:4}}]=0.24}$ to scale the measurements. The number of viable cells is calculated by inserting the estimated expected values of the model parameters~$\prol$, $\Kv$, $m$ and~$\nat$ from Table~\ref{Tab:deviations1} into the analytical solution~\eqref{eq:solV}. Figure~\ref{Fig:resOpt} compares the corresponding scaled data with the calculated time course of~$V(t)$.
As expected, both solutions are very similar, as for model~\ref{eq:ODEprol} only the parameters~$\prol$, $\Kv$, $m$, and~$\nat$ are relevant, which have similar estimates for both~\ref{eq:ODE} and~\ref{eq:ODEstress}. Furthermore, we see a good fit to the data (see also ratios of validation metrics in Table~\ref{Tab:validResSUP} in the supplement). It is observable that the model solutions tend to overestimate the data after the steady state is reached. In reality, population growth often follows the concept of ``overshoot'': this describes the observation that a population might exceed its carrying capacity and drop back down afterwards (see e.g.~\cite{overshoot}). Reasons can be a lag in sensing the available resources in the environment or the time interval between birth and death. The measurements starting with~${\initV\in\left\{1.00,\,0.50,\,0.25\right\}}$ (orange, purple, and green markers in Figure~\ref{Fig:resOpt}) suggest an observable overshoot effect, but the model solutions cannot reproduce this phenomenon due to the limited structure of model~\ref{eq:ODEprol}. This could explain the discrepancy between solution and data close to the carrying capacity.
\section{Conclusion}
\label{sec:conclusion}
In summary, the presented results show that model~\ref{eq:ODEstress} is an appropriate alternative way to describe the influences of nutrient changes on viable cells. The model preserves positivity and boundedness of the solutions and calibration yields biologically reasonable parameter values. Keeping in mind that large deviations from the data sets~D1--D5 focus on a particular measurement (${t\approx 3}$), we altogether see a good agreement to the data for different nutrient scenarios. For the simple setting with the nutrients as the only environmental variable, there is no strong need from a quantitative perspective to model its influence by using the additional variable~$\stress$. But even in this  special case, the Bayes factor indicates a preference for including the environmental stress level $\eta$ into modeling. A mathematical model should be kept as simple as possible, but as additional influencing factors may become relevant, a suitable and flexible modeling interface is important. From a qualitative perspective, the general modeling approach~\eqref{eq:stress} can provide such a setup.

The introduction of an artificial stress level to describe different environmental factors in a collective manner has several advantages. Due to the modular structure of the general ODE for the stress level~$\stress$ (equation~\ref{eq:stress}), this model is easily extendable to several environmental variables. This enables the consideration of different sensitivities of the cells to changes in specific environmental factors. If the values of the corresponding sensitivity rates~$\alpha^\pm_j$ can be recovered from data, the magnitude of them can give insight into time scales of the cells' reaction to different environmental changes. This information can be used on the one hand to simplify the model by quasi-steady state assumptions for environmental variables regarding large sensitivity rates. On the other hand, if some sensitivity rates are estimated to be close to zero, the model can be reduced easily by setting them to zero and neglecting the ODEs of the corresponding environmental variables. The resulting quality of fit of the simplified model can then be investigated e.g. by the Bayes factor.

Overall, the model allows to describe a survival strategy of the cells assuming they adjust their behavior to avoid a stressful environment by prioritizing according to their sensitivities. If spatial inhomogeneity is considered (i.e. if the model is extended to a system of partial differential equations), such an avoidance strategy can be modeled by influencing the cells' motility. This can be mathematically described similarly to chemotaxis: cells might try to move preferably to areas with a low stress level to improve their growth/survival conditions. Future work aims to investigate the application of the proposed general model in such a setting, especially considering mechanical environmental features like pressure or stiffness of the surrounding tissue, as considered experimentally e.g. in \cite{Ozkan1}. New experimental approaches like ``organ-on-a-chip'' will allow for even more realistic and nevertheless (as {\it in vitro} technique) well-controlled settings and corresponding parameter calibration \cite{Ozkan2, Ozkan3}.
\section*{Acknowledgments}
This paper was supported by Deutsche Forschungsgemeinschaft (DFG) through TUM International Graduate School of Science and Engineering (IGSSE), GSC 81.

\section*{Conflict of interest}
All authors declare no conflicts of interest in this paper.

\beginsupplement
\section{Supplementary}
\subsection{Analytical solution of model~\ref{eq:ODE}}
\label{SUPsec:analSol}
The following equations~\eqref{eq:analSol1} and~\eqref{eq:analSol3} are the case-dependent analytical solutions of system~\eqref{eq:ODE} as given in Section~\ref{subsec:model}. The corresponding calculations are done in the following paragraphs. It can be shown by differentiation that these derived terms for~$V(t)$ actually solve the given initial value problem.
\paragraph{Case 1: $\boldsymbol{\scaledProl>\scaledDeath}$.} Model~\ref{eq:ODE} has the same form as model~\ref{eq:ODEprol} just with nutrient-scaled growth rate~${\scaledProl=\actFct{\initS}\cdot\prol}$ and death rate~${\scaledDeath=\nat+\deactFct{\initS}\cdot\starv}$\,, each depending on the constant nutrient supply~$\initS$. Hence, for ${\scaledProl>\scaledDeath}$ the ODE of~\eqref{eq:ODE} can again be written in the form of a logistic growth term and the analytical solution can be given analogously to equation~\eqref{eq:solV} with scaled rates:	\begin{equation}
	V(t)=\Bigg(\,\frac{(\initV\Kv)^m(\scaledProl-\scaledDeath)}{{\scaledProl}\initV^m+\left(\left(\scaledProl-\scaledDeath\right)\Kv^m-\scaledProl\initV^m\right)\cdot e^{-\left(\scaledProl-\scaledDeath\right)mt}}\,\Bigg)^{\frac{1}{m}}\,. \label{eq:analSol1}
\end{equation}
\paragraph{Case 2: $\boldsymbol{\scaledProl<\scaledDeath}$.} In case of~${\scaledProl\leq\scaledDeath}$ the ODE cannot be rewritten in pure logistic form, since the root's argument of the resulting ``net'' carrying capacity~${\Kv(1-\scaledDeath/\scaledProl)^{\frac{1}{m}}}$ might cause problems. However, it can be shown that for~${\scaledProl<\scaledDeath}$ the analytical solution above is still well-defined, since the root's argument is non-negative:
\begin{equation}
	V(t)\overset{\eqref{eq:analSol1}}{=}\Bigg(\,\frac{(-1)}{(-1)}\cdot\frac{\overbrace{(\initV\Kv)^m}^{\geq 0}\overbrace{(\scaledDeath-\scaledProl)}^{>\,0}}{\underbrace{\scaledProl\initV^m}_{\geq\,0}\cdot\big(-1+\underbrace{e^{(\scaledDeath-\scaledProl)mt}}_{>\,1}\big)+\underbrace{(\scaledDeath-\scaledProl)\Kv^m}_{>\,0}\cdot \underbrace{e^{(\scaledDeath-\scaledProl)mt}}_{>\,1}}\,\Bigg)^{\frac{1}{m}}\,.\label{eq:analSol2}
\end{equation} 
\paragraph{Case 3: $\boldsymbol{\scaledProl=\scaledDeath}$.} For ${\scaledDeath\to\scaledProl}$ the fraction inside the root of formula~\eqref{eq:analSol2} degenerates, as numerator and denominator both tend to zero. Since $V(t)$ is a differentiable function in the parameter space, we can derive the analytical solution in this case using l'Hospital's rule (l'H). It holds
\begin{align}
	\lim\limits_{\scaledDeath\to\scaledProl}& V(t)^m\overset{\hphantom{\text{l'H}}}{=}\lim\limits_{\scaledDeath\to\scaledProl}~
	\frac{(\initV\Kv)^m\overbrace{(\scaledDeath-\scaledProl)}^{\to 0}}{\scaledProl\initV^m\cdot\underbrace{\big(e^{(\scaledDeath-\scaledProl)mt}-1\big)}_{\to 0}+\underbrace{(\scaledDeath-\scaledProl)}_{\to 0}\Kv^m\cdot \underbrace{e^{(\scaledDeath-\scaledProl)mt}}_{\to1}}\nonumber \\
	&\overset{\text{l'H}}{=}\lim\limits_{\scaledDeath\to\scaledProl}~
	\frac{(\initV\Kv)^m}{\scaledProl\initV^m\cdot mt\cdot e^{(\scaledDeath-\scaledProl)mt}+\Kv^m\cdot e^{(\scaledDeath-\scaledProl)mt}+(\scaledDeath-\scaledProl)\Kv^m\cdot mt\cdot e^{(\scaledDeath-\scaledProl)mt}} \nonumber\\
	&\overset{\hphantom{\text{l'H}}}{=} \lim\limits_{\scaledDeath\to\scaledProl}~
	\frac{(\initV\Kv)^m}{\Big(mt\cdot\big(\scaledProl\initV^m+\underbrace{(\scaledDeath-\scaledProl)}_{\to 0}\Kv^m\big)+\Kv^m\Big)\cdot \underbrace{e^{(\scaledDeath-\scaledProl)mt}}_{\to 1}}=\frac{(\initV\Kv)^m}{mt\scaledProl\initV^m+\Kv^m}\,,\nonumber \\
	\mbox{yielding}~&\lim\limits_{\scaledDeath\to\scaledProl}V(t)\overset{\hphantom{\text{l'H}}}{=}\initV\Kv\cdot\big(mt\scaledProl\initV^m+\Kv^m\big)^{-1/m}\,. \label{eq:analSol3}
\end{align} 
\subsection{Mathematical analysis of the ODE models}
\label{SUPsec:mathAna}
In the following paragraphs, the solutions of all models are analyzed in terms of positivity and boundedness. Furthermore, we give the respective steady states and check for their stability.
\paragraph{Positivity and boundedness of the solutions.} In a reasonable biological context, the size of a cell population cannot grow arbitrarily because of environmental limitations. Hence, for each of the presented models to be feasible, the corresponding solution $V(t)$ has to stay bounded and especially non-negative. The same properties are desired for the environmental stress level~$\stress$.

\smallskip

\textsl{Environmental stress level}. Before analyzing the variable~$V$ in all models, we start with the stress level~$\stress$ from system~\eqref{eq:ODEstress}. Given ${\stress(0)=\initStress\in[0,1]}$, its ODE yields
 \begin{equation}
	\dStress=\adaptS\big(\deactFct{\initS}-\stress\big)~\begin{cases}
		~ > 0 & \mbox{for } \stress< \deactFct{\initS}\\
		~ = 0 & \mbox{for } \stress= \deactFct{\initS}\\
		~ < 0 & \mbox{for } \stress> \deactFct{\initS}
	\end{cases}\label{eq:stressLvl}\,,
\end{equation} 
leading to both positivity and boundedness of~$\stress$: ${0\leq\stress\leq\max\left\{\initStress~,~\deactFct{\initS}\right\}\leq 1}$. Analogous bounds can be calculated for the general ODE~\eqref{eq:stress} for the environmental stress level, since ${\alpha^\pm_i,\,\delta^\pm_j\geq 0}$.

\smallskip

\textsl{Density of viable cells}. Regarding the remaining variable~$V$, positivity and boundedness in model~\ref{eq:ODEprol} follow directly from the fact that its ODE can be rewritten as a logistic growth model with positive growth rate~$\netProl$ and capacity~$\netKv$, preserving positivity and boundedness by definition. To conclude these properties also for models~\ref{eq:ODE} and~\ref{eq:ODEstress}, we analyze a more general, non-autonomous model of the form
\begin{align}
	&\!\left\{~\begin{aligned}
		\dV &= \underbrace{h(t)\cdot\prol}_{=\,\prol_h(t)}\, V\left(1-\left(\frac{V}{\Kv}\right)^{m}\right)-\underbrace{\Big(\nat+\big(1-h(t)\big)\cdot\starv\Big)}_{=\,\nat_h(t)} V\,, \\
		V(0)&=\initV\,,
	\end{aligned} \right. & \label{eq:varS}
\end{align} 
where $h$ is a potentially time-dependent and bounded function:  ${h(t)\in[0,1]~\forall t\geq 0}$. In fact, problem~\eqref{eq:varS} is a generalization of both systems~\eqref{eq:ODE} and~\eqref{eq:ODEstress}, where respectively ${h(t)=\actFct{\initS}}$ and ${h(t)=1-\stress(t)\overset{\eqref{eq:stressSol}}{=}1-\deactFct{\initS}\cdot(1-e^{-\adaptS t})-\initStress e^{-\adaptS t}}$. For ${h(t)=\actFct{\initS}}$ the required boundedness of~$h$ is given by definition of ${\delta^+(S)}$ and for ${h(t)=1-\stress(t)}$ it is given by the previously derived bounds: ${\stress\in[0,1]}$. Hence, we can analyze the general problem~\eqref{eq:varS} for positivity and boundedness instead and transfer the results to models~\ref{eq:ODE} and~\ref{eq:ODEstress}.

With ${h(t)\in[0,1]~\forall t\geq 0}$ and ${\initV\geq 0}$, positivity follows directly by observing ${\dV =0\geq 0}$ for ${V=0}$. To analyze boundedness, we check the sign of the derivative of~$V$. For ${\nat_h(t)\geq \prol_h(t)}$ we estimate:
\begin{equation*}
\dV\overset{\eqref{eq:varS}}{=}\prol_h(t)\, V\left(1-\left(\frac{V}{\Kv}\right)^{m}\right)-\nat_h(t) V
=\prol_h(t)\,V\,\Bigg(\underbrace{1-\frac{\nat_h(t)}{\prol_h(t)}}_{\leq 0}-\left(\frac{V}{\Kv}\right)^m\Bigg)\leq0\,.
\end{equation*} 
In case of~${\nat_h(t)<\prol_h(t)}$, we can rewrite the latter term further by placing $\left(1-\frac{\nat_h(t)}{\prol_h(t)}\right)$ outside the bracket:
\begin{equation}
\dV\!=\underbrace{\big(\prol_h(t)-\nat_h(t)\big)}_{>0}V\Bigg(\!1-\Bigg(\frac{V}{\Kv\underbrace{\big(1-\nat_h(t)/\prol_h(t)\big)^{\frac{1}{m}}}_{>0}}\,\Bigg)^{\!\!m}\Bigg)\begin{cases}
	 > 0 & \mbox{if } V< \Kv\left(1-\frac{\nat_h(t)}{\prol_h(t)}\right)^{\frac{1}{m}},\\
	 = 0 & \mbox{if } V= \Kv\left(1-\frac{\nat_h(t)}{\prol_h(t)}\right)^{\frac{1}{m}},\\
	 < 0 & \mbox{if } V>\Kv\left(1-\frac{\nat_h(t)}{\prol_h(t)}\right)^{\frac{1}{m}}.
\end{cases}\label{eq:capacity}
\end{equation} 
Together with ${\prol_h(t)=h(t)\cdot\prol\leq\prol}$ and~${\nat_h(t)=\nat+\big(1-h(t)\big)\cdot\starv\geq\nat}$\,, this eventually leads to a time-independent upper bound for the population size:
\begin{equation*}
	V\leq\max\left\{\initV~,~ \Kv\left(1-\frac{\nat_h(t)}{\prol_h(t)}\right)^{\frac{1}{m}}\right\}\leq\max\left\{\initV~,~ \Kv\left(1-\frac{\nat}{\prol}\right)^{\frac{1}{m}}\right\}\,.
\end{equation*} 
Note that relations~\eqref{eq:stressLvl} and~\eqref{eq:capacity} also fit to the non-trivial steady states of the corresponding models (see next paragraph).
\paragraph{Steady states of the models.} To calculate the steady states of each problem, we introduce the following notations for the right hand sides of the corresponding autonomous ODEs in each model:
\begin{equation}
\eqref{eq:ODEprol}:~ \dV=f_{\text{opt}}(V)\,,\quad
\eqref{eq:ODE}:~ \dV=f_S(V)\,,\quad
\eqref{eq:ODEstress}:~ \dV=f_\stress(V,\stress)\,,~ \dStress = g_{\stress}(V,\stress)\,.\label{eq:RHS}
\end{equation} 
The steady states are derived by setting the right hand sides to zero respectively. From a biological point of view, we would expect at least two steady states~$\equibV$ for the cell density. The cells can either go extinct or the population size tends to the ``net'' carrying capacity, which depends on the relation between growth and death rate. For a growing population, we use more general notations during the analysis: Let $b$ be the total growth rate, which is larger than the total death rate~$a$ and define the function
\begin{equation*}
	\mathcal K:\left\{\,(a,b)\in\mathbb{R}_+^2~\big|~a<b\,\right\}\to(0,\Kv) \quad\mbox{with}\quad (a,b)\mapsto\Kv\big(1-a/b\big)^{\frac{1}{m}}
\end{equation*} 
to be the ``net'' carrying capacity depending on the rates. Then, the mentioned expected steady states translate to ${\equibV=0}$ or ${\equibV=\mathcal K(a,b)}$.
\smallskip

\textsl{System~\eqref{eq:ODEprol}}. For model~\ref{eq:ODEprol} these are actually the only two steady states:
\begin{align*}
	\dV|_{V=\equibV}\overset{\eqref{eq:ODEprol}}{=}&\prol\,\equibV \left(1-\left(\frac{\equibV }{\Kv}\right)^m\right)-\nat\equibV =\equibV\cdot \left(\prol \left(1-\left(\frac{\equibV }{\Kv}\right)^m\right)-\nat\right)=0 \label{eq:equibV}\\[2pt]
	&\Rightarrow\quad \equibV =0\quad\mbox{or}\quad~ \prol\left(1-\left(\frac{\equibV }{\Kv}\right)^m\right)-\nat = 0 \quad\overset{\prol>\nat}{\Longrightarrow}~ \equibV  = \Kv\bigg(1-\frac{\nat}{\prol}\bigg)^{\frac{1}{m}}=\mathcal K(\nat,\prol)\,.
\end{align*} 
Note that the same non-trivial steady state can also be calculated by considering the limit~${t\to\infty}$ of the corresponding analytical solution~\eqref{eq:solV}.
\smallskip

\textsl{System~\eqref{eq:ODE}}. Since model~\ref{eq:ODE} only differs from model~\ref{eq:ODEprol} by time-independent scaling of the growth and death rate, analogous calculations lead to a similar result:
\begin{equation*}
	\dV|_{V=\equibV}\!\overset{\eqref{eq:ODE}}{=}\!\scaledProl\equibV\!\left(\!1-\left(\frac{\equibV }{\Kv}\right)^{\!\!m}\right)-\scaledDeath\equibV\!=0~\Rightarrow\,\equibV\!=\begin{cases}
		0 & \mbox{if }\scaledDeath\geq\scaledProl, \\
		0\mbox{ or }\Kv\left(1-\frac{\scaledDeath}{\scaledProl}\right)^{\frac{1}{m}}\!=\mathcal K(\scaledDeath,\scaledProl) & \mbox{if } \scaledDeath<\scaledProl.
	\end{cases}
\end{equation*} 
Hence, we have the two expected steady states in case of~${\scaledDeath<\scaledProl}$, but only the trivial one if~${\scaledDeath\geq\scaledProl}$.
\smallskip

\textsl{System~\eqref{eq:ODEstress}}. For the last model~\ref{eq:ODEstress}, we have to consider a system of two ODEs:
\begin{align*}
	&\!\eqref{eq:ODEstress}:\quad\left\{~\begin{aligned}
		\dV|_{V=\equibV,\,\stress=\bar{\stress}} &= (1-\bar{\stress})\cdot\prol\, \equibV\left(1-\left(\frac{\equibV}{\Kv}\right)^{m}\right)-(\nat+\bar{\stress}\cdot\starv) \equibV=0\,, \\
		\dStress|_{V=\equibV,\,\stress=\bar{\stress}} &= \adaptS\deactFct{\initS}\big(1-\bar{\stress}\big)-\adaptS\bar{\stress}\actFct{\initS}=\adaptS\big(\deactFct{\initS}-\bar{\stress}\big)=0\,.
	\end{aligned} \right. & 
\end{align*}
Latter equation is independent of~$V$, which directly results in the steady state for the stress level: ${\bar{\stress}=\deactFct{\initS}}$. Inserting this into the first equation makes it identical to the equation from the previous paragraph for system~\eqref{eq:ODE}. Therefore, the steady states of problem~\eqref{eq:ODEstress} are
\begin{align*}
\left(\begin{array}{c}
	\equibV \\ \bar{\stress}	
	\end{array}\right)&=\left(\begin{array}{c}
	0 \\ \deactFct{\initS}	
\end{array}\right)\mbox{ for all }\scaledDeath,\scaledProl\\
\mbox{and}\quad\left(\begin{array}{c}
\equibV \\ \bar{\stress}	
\end{array}\right)
&=\left(\begin{array}{c}
\Kv\left(1-\frac{\scaledDeath
}{\scaledProl}\right)^{\frac{1}{m}} \\ \deactFct{\initS}	
\end{array}\right)=\left(\begin{array}{c}
\mathcal K(\scaledDeath,\scaledProl) \\ \deactFct{\initS}	
\end{array}\right),\mbox{ if }\scaledDeath<\scaledProl\,.
\end{align*}
\paragraph{Stability of the steady states.} The derived steady states can be analyzed for stability. We use the previously introduced notations~\eqref{eq:RHS} for the right hand sides of the ODEs. Following the theorem of Hartman and Grobman~\cite{Hartman, Grobman}, local stability of a steady state can be evaluated by checking the sign of the derivatives of the right hand sides in these points.
\smallskip

\textsl{Systems~\eqref{eq:ODEprol} and~\eqref{eq:ODE}}. For the first two models we have to check the sign of $f'_{\text{opt}}(\equibV )$ and $f'_S(\equibV )$ respectively. Keeping the constraint ${\prol>\nat}$ in mind, inserting the steady states of model~\ref{eq:ODEprol} yields
\begin{equation*}
	f_{\text{opt}}(V)=\odeVprol\quad\Rightarrow\quad f'_{\text{opt}}(V)=\prol-\nat-\frac{(m+1)\prol}{\Kv^m}V^m
\end{equation*}
resulting in
\begin{equation*}
	f'_{\text{opt}}(0)=\prol-\nat>0\quad
	\mbox{and}\quad f'_{\text{opt}}\left(\Kv\bigg(1-\frac{\nat}{\prol}\bigg)^{\frac{1}{m}}\right)=\prol-\nat-(m+1)(\prol-\nat)=-m(\prol-\nat)<0\,.
\end{equation*} 
Hence, the non-trivial steady state is stable, whereas the trivial one is not. Since $f_S$ has the same form as $f_{\text{opt}}$ only with scaled growth and death rate, analogous calculations result in
\begin{equation*}
	f'_S(0)=
	\scaledProl-\scaledDeath\begin{cases}
	> 0 &\mbox{if }\scaledDeath<\scaledProl\,, \\ 
	= 0 &\mbox{if }\scaledDeath=\scaledProl\,, \\ 
	< 0 &\mbox{if }\scaledDeath>\scaledProl\,. 
	\end{cases}~
	\mbox{and, for }\scaledDeath<\scaledProl\!:~ f'_S\left(\!\Kv\bigg(1-\frac{\scaledDeath}{\scaledProl}\bigg)^{\!\!\frac{1}{m}}\right)\!\!=-m(\scaledProl-\scaledDeath)< 0\,.
\end{equation*} 
The stability of the steady states now depends on the relation between the scaled growth rate~$\scaledProl$ and death rate~$\scaledDeath$. For the first case~${\scaledDeath<\scaledProl}$, the trivial steady state is locally unstable, whereas the non-trivial one is stable. For~${\scaledDeath>\scaledProl}$ only the trivial steady state exists, which is stable.

In case of equality ${\scaledDeath=\scaledProl}$, the theorem of Hartman and Grobman cannot be used to determine stability for~${\equibV=0}$ and we take another approach by investigating the change of a time-dependent perturbation from the steady state~${\equibV =0}$. Observing that ${f_S(V)}$ is a polynomial in~$V$ of degree ${m+1}$ without an absolute term, all but one derivatives vanish in~${\equibV =0}$:
\begin{equation}
	f_S^{(j)}(0)=~\begin{cases}
	~-\frac{\scaledProl}{\Kv^m}\cdot\frac{(m+1)!}{(m+1-j)!}\cdot 0^{m+1-j}=0\qquad &\mbox{for }j<m+1\,, \\[2pt]
	~-\frac{\scaledProl}{\Kv^m}(m+1)!=\text{const.}\neq 0 &\mbox{for } j=m+1\,, \\[2pt]
	~0 &\mbox{for } j>m+1\,.
\end{cases} \label{eq:derivRHS}
\end{equation}
Let~${\nu=\nu(t)>0}$ be the time-dependent perturbation from the steady state~${\equibV =0}$ with initial value~${\nu(0)=\nu_0}$. The change of the perturbation in time can be derived by insertion into the ODE, i.e. $f_S\big(\equibV +\nu(t)\big)\overset{\text{ODE}}{=}\totDiff{t}\big(\equibV +\nu(t)\big)=\totDiff{t}\nu$, and Taylor series expansion of the right hand side:
\begin{align*}
\totDiff{t}\nu&\overset{\hphantom{(0.0)}}{=}f_S\big(\equibV +\nu(t)\big)\\
&\overset{\eqref{eq:derivRHS}}{=}\left(\sum_{j=0}^{m} \underbrace{f_S^{(j)}(0)}_{=\,0} \frac{(\nu-0)^j}{j!}\right)+\left(\frac{f_S^{(m+1)}(0)}{(m+1)!}(\nu-0)^{m+1}\right)+\left(\sum_{j=m+2}^{\infty}\underbrace{f_S^{(j)}(0)}_{=\,0}\frac{(\nu-0)^j}{j!}\right) \\&\overset{\hphantom{(0.0)}}{=} -\frac{\scaledProl}{\Kv^m}\nu^{m+1}.
\end{align*} 
The derived ODE in~$\nu(t)$ can be analytically solved by separation of variables, resulting in
\begin{equation*}
	\nu(t)=\nu_0\Kv\cdot\left(mt\scaledProl\nu_0^m+\Kv^m\right)^{-1/m}\,.
\end{equation*} 
Taking the limit~${t\to\infty}$ shows that the perturbation vanishes over time, i.e. local stability of the steady state~${\equibV =0}$ can be concluded.
\smallskip

\textsl{System~\eqref{eq:ODEstress}}. For the last model~\ref{eq:ODEstress} we again use Hartman-Grobman theorem to investigate the corresponding locally linearized model. Reminding that the system is given by
\begin{align*}
	&\!\eqref{eq:ODEstress}:\quad\left\{~\begin{aligned}
		\dV &= f_\stress(V,\stress) =(1-\stress)\cdot\prol\, V\left(1-\left(\frac{V}{\Kv}\right)^{m}\right)-(\nat+\stress\cdot\starv) V\,, \\
		\dStress &=g_\stress(V,\stress)= \adaptS\deactFct{\initS}(1-\stress)-\adaptS\actFct{\initS}\stress=\adaptS\big(\deactFct{\initS}-\stress\big)\,,
	\end{aligned} \right. & 
\end{align*} 
linearization in the steady state~$(\equibV,\bar{\stress})$ with the Jacobian~$J$ yields
\begin{align*}
	\left(\begin{array}{c}
		\dot{x} \\ \dot{y}
	\end{array}\right)
	&=J\left(\begin{array}{c}
		{x} \\ {y}
	\end{array}\right)=\left.\left(\begin{array}{cc}
	\partDiff{V}\,f_\stress & \partDiff{\stress}\, f_\stress \\
	\partDiff{V}\, g_\stress & \partDiff{\stress}\, g_\stress
\end{array}\right)\right|_{(V,\stress)=(\equibV,\bar{\stress})}\left(\begin{array}{c}
{x} \\ {y}
\end{array}\right)\\
	&=\left(\begin{array}{cc}
(1-\bar{\stress})\cdot\prol\left(1-\frac{m+1}{\Kv^m}\equibV^m\right)-(\nat+\bar{\stress}\cdot\starv) & -\prol\,\equibV\left(1-\left(\frac{\equibV}{\Kv}\right)^m\right)-\starv\,\equibV \\
0 & -\adaptS
\end{array}\right)\left(\begin{array}{c}
{x} \\ {y}
\end{array}\right)\,.
\end{align*} 
The eigenvalues~$\chi_{1,2}$ of the triangular Jacobian~$J$ in dependence on the steady state~$(\equibV,\bar{\stress})$ are given by
\begin{equation*}
	\chi_1=-\adaptS<0\qquad\mbox{and}\qquad\chi_2=\chi_2(\equibV,\bar{\stress})=(1-\bar{\stress})\cdot\prol\left(1-\frac{m+1}{\Kv^m}\equibV^m\right)-(\nat+\bar{\stress}\cdot\starv)\,.
\end{equation*} 
The first eigenvalue~$\chi_1$ is always negative, whereas the sign of~$\chi_2$ depends on the relation between~$\scaledProl$ and~$\scaledDeath$. In particular, the non-trivial steady state only exists for ${\scaledDeath<\scaledProl}$, leading to
\begin{equation*}
	\chi_2\left(\Kv\left(1-\frac{\scaledDeath}{\scaledProl}\right)^{\frac{1}{m}},\deactFct{\initS}\right)
	=\scaledProl\left(1-(m+1)\left(1-\frac{\scaledDeath}{\scaledProl}\right)\right)-\scaledDeath=-m(\scaledProl-\scaledDeath)<0\,,
\end{equation*} 
i.e. this steady state is a stable node. For the trivial steady state it holds
\begin{equation*}
	\chi_2\big(0,\deactFct{\initS}\big)
	=\scaledProl-\scaledDeath\,\begin{cases}
	\,>0 &\mbox{for }\scaledDeath<\scaledProl\,, \\[-3pt]
	\,=0 &\mbox{for }\scaledDeath=\scaledProl\,, \\[-3pt]
	\,<0 &\mbox{for }\scaledDeath>\scaledProl\,, \\
\end{cases}
\end{equation*}
and, therefore, this steady state is a stable node for ${\scaledDeath>\scaledProl}$ and an unstable one for ${\scaledDeath<\scaledProl}$. In case of equality~${\scaledDeath=\scaledProl}$ the second eigenvalue gets zero and the theorem is again not applicable anymore. To show stability in this case, we use the fact that the ODE of $\stress$ is independent of $V$, i.e. ${g_\eta(V,\eta)=g_\eta(\eta)}$, and thus the equations of model~\ref{eq:ODEstress} can be analyzed consecutively. On the one hand, it holds
\begin{equation*}
	g_\stress(\stress)=\adaptS\big(\deactFct{\initS}-\stress\big)\quad\Rightarrow \quad g'_\stress(\stress)=-\adaptS<0
\end{equation*} 
and on the other hand, inserting ${\bar{\stress}=\deactFct{\initS}}$ into the first ODE of the system yields
\begin{equation*}
	f_\stress(V,\bar\stress)=\scaledProl\,V \big(1-(V/\Kv)^m\big)-\scaledDeath\,V=f_S(V)\,.
\end{equation*} 
For this equation stability was already shown for model~\ref{eq:ODE} in case of~${\scaledDeath=\scaledProl}$. Eventually, the trivial state of model~\ref{eq:ODEstress} is stable for this case as well.
\subsection{Further model calibration results}
\label{SUPsec:calRes}
\begin{figure}[H]
	\begin{center}
		\includegraphics[width=.8\textwidth]{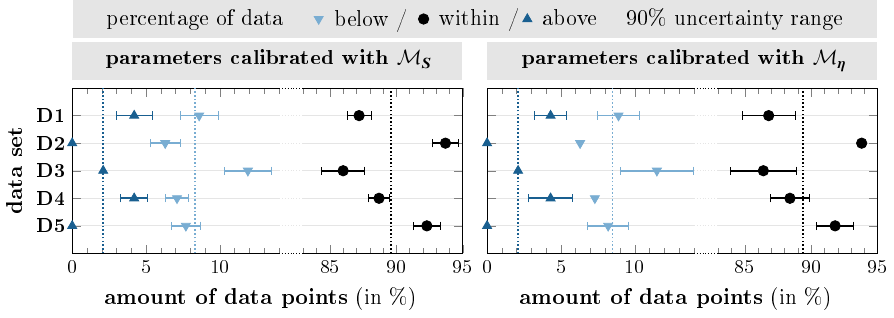}
		\caption{Percentage of data points which are situated below/within/above the $90\%$ range ${[V\cdot P_{5\%}\,,\, V\cdot P_{95\%}]}$ around the solution~$V$. We remind that the error bars indicate the numerical deviation between the SMC runs. The dotted lines in the corresponding color give the average over all data sets~D1--D5}
		\label{Fig:resNoise90}
	\end{center}
\end{figure}

\paragraph{Calibrated parameters.} For the proportionality constants and reparametrization parameters, the expected values and variances of the posteriors are listed in Table~\ref{Tab:deviationsSUP}.
\begin{table}[H]
	\setlength{\tabcolsep}{4pt}
	\begin{center}
		\caption{Mean values and variances of the posterior distributions for some calibrated parameters. The deviations between the runs of the SMC algorithm are given in terms of the $95\%$ prediction interval of a normal distribution~$\mathcal{N}(\mu,\sigma^2)$: ${\mu\pm 1.96\sigma}$.}
		\begin{tabular}{c r c c c c c}
			\hline
			~ & ~ & ~ & ~ & ~ & \llap{\parbox{8.65cm}{\textbf{Calibrated parameters}}} \\[-2pt]
			\textbf{Model} & ~ & $\boldsymbol{n_{\text{\textbf{D1:4}}}}$ & $\boldsymbol{n_{\text{\textbf{D5}}}}$ & $\boldsymbol{c_n}$ & $\boldsymbol{c_1}$ & $\boldsymbol{c_2}$ \\
			\hline
			\ref{eq:ODE} & $\boldsymbol{\mean(\,\cdot\,)}$ & $0.244\pm0.006$ & $0.190\pm0.020$ & $0.780\pm0.084$ & $0.236\pm0.076$ & $0.576\pm0.293$ \\
			\ref{eq:ODEstress} & ~ & $0.243\pm0.004$ & $0.182\pm0.022$ & $0.752\pm0.094$ & $0.240\pm0.046$ & $0.554\pm0.159$ \\[4pt]
			\ref{eq:ODE} & \textbf{Var}$\boldsymbol{(\,\cdot\,)}$ & $0.002\pm0.004$ & $0.008\pm0.017$ & $0.034\pm0.069$ & $0.011\pm0.024$ & $0.053\pm0.100$ \\
			\ref{eq:ODEstress} & ~ & $0.002\pm0.004$ & $0.007\pm0.014$ & $0.032\pm0.055$ & $0.016\pm0.036$ & $0.069\pm0.122$  \\
			\hline
		\end{tabular}
		\label{Tab:deviationsSUP}
	\end{center}
\end{table}
\paragraph{Uncertainties in the model.} In Figure~\ref{Fig:resNoise90} we consider the 90\% uncertainty range~\eqref{eq:range} with the percentiles~$P_{5\%}$ and~$P_{95\%}$ calculated in Section~\ref{subsec:resComp}, using either the calibrated parameters of \ref{eq:ODE} or \ref{eq:ODEstress} (left/right plot). There are three different markers per data set (vertical axis), showing how many data points (horizontal axis) are situated below/within/above the $90\%$~range. For instance regarding D1 (topmost on vertical axis), in both plots the black bullet marker shows that approx. 87\% of the measurements actually lie within the 90\% range, whereas the dark/light blue triangle marker shows that about 4\%/9\% of them are above/below the 90\% range. 
\paragraph{Quality of fit to the data.} The validation metric~\eqref{eq:valid} measures the mismatch between data and model prediction, therefore
\begin{equation*}
	\frac{d_\stress}{d_S}=\frac{d\Big(F^{\text{data}},F^{\text{sol}}_{\mathcal{M}_\stress}\Big)}{d\Big(F^{\text{data}},F^{\text{sol}}_{\mathcal{M}_S}\Big)}~\begin{cases}
		~\ll 1 & \text{indicates: better fit for~\ref{eq:ODEstress}}\,, \\
		~\approxeq 1 & \text{indicates: \ref{eq:ODE} and~\ref{eq:ODEstress} fit equally well}\,, \\
		~\gg 1 & \text{indicates: better fit for~\ref{eq:ODE}}\,.
	\end{cases}
\end{equation*} 
Calculations for each data set~D1--D6 and initial cell number~$\initV$ result in the values given in Table~\ref{Tab:validResSUP}. For data set~D6 the solution of model~\ref{eq:ODEprol} was used by inserting the average estimated parameters resulting from the calibration of~\ref{eq:ODE} resp.~\ref{eq:ODEstress}.
\begin{table}[H]
	\setlength{\tabcolsep}{4pt}
	\begin{center}
		\caption{Ratio~$d_\stress/d_S$ of the validation metrics using model~\ref{eq:ODEstress} resp.~\ref{eq:ODE} averaged over all time points considering measurements regarding various $\initV$ (columns) and data sets, i.e. different $\initS$ (rows). The values in the last row/column are averaged over all data sets/initial cell densities before taking the mean over the SMC runs. The numerical deviations between the runs are given in terms of the $95\%$ confidence interval of a normal distribution.}
		\begin{tabular}{ccccccc} \hline
			~ & & & & & & \llap{\parbox{12.9cm}{\textbf{Initial cell density }$\boldsymbol{\initV}$}} \\[-2pt]
			\textbf{Data set} & \textbf{1.00}&\textbf{0.50}&\textbf{0.25}&\textbf{0.10}&\textbf{0.05} & \quad \textbf{all (avg.)}\\\hline
			\textbf{D6} & $1.05\pm0.54$ & $1.04\pm0.33$& $1.07\pm0.17$ & $1.01\pm0.24$ & $1.03\pm0.13$ &\quad$1.04\pm0.26$\\[5pt]
			\textbf{D1} & $1.04\pm0.14$ & $1.01\pm0.12$& $1.00\pm0.09$ & - & - &\quad$1.02\pm0.04$\\
			\textbf{D2} & $1.01\pm0.07$ & $0.99\pm0.17$& $0.99\pm0.15$ & - & - &\quad$1.00\pm0.12$\\
			\textbf{D3} & $1.02\pm0.04$ & $1.02\pm0.21$& $1.03\pm0.13$ & - & - &\quad$1.02\pm0.07$\\
			\textbf{D4} & $1.00\pm0.06$ & $1.01\pm0.13$& $1.01\pm0.04$ & - & - &\quad$1.01\pm0.07$\\
			\textbf{D5} & $0.94\pm0.15$ & $1.06\pm0.18$& $1.03\pm0.08$ & - & - &\quad$1.01\pm0.11$\\[5pt]
			\textbf{all (avg.)} & $1.01\pm0.09$ & $1.02\pm0.12$& $1.02\pm0.04$ & $1.01\pm0.24$ & $1.03\pm0.13$ &\quad$1.02\pm0.09$\\
			\hline
		\end{tabular}
		\label{Tab:validResSUP}
	\end{center}
\end{table}
All values are  close to one and the differences between them are in the scale of numerical variations of the SMC algorithm. Hence, the investigation of the validation metric does not show any preference of a particular model.
\par
\end{document}